\documentclass[11pt]{article}  
\usepackage{amsmath}
\usepackage{amssymb}
\usepackage{theorem}
\usepackage{pstricks}
\usepackage{xcolor}
\title{\sffamily Functions with Prescribed Best Linear Approximations}
\author{Patrick L. Combettes$^1$
and Noli N. Reyes$^2$ \\[5mm]
$\!^1$UPMC Universit\'e Paris 06\\
Laboratoire Jacques-Louis Lions -- UMR 7598\\
75005 Paris, France\\
{\ttfamily plc@math.jussieu.fr}\\[4mm]
$^2$University of the Philippines -- Diliman\\
Institute of Mathematics\\
Quezon City, 1101 Philippines\\
{\ttfamily noli@math.upd.edu.ph}}
\date{~}
\newcommand{\scal}[2]{{\langle{{#1}\mid{#2}}\rangle}}

\newcommand{\abscal}[2]{\left|\left\langle{{#1}\mid{#2}}%
\right\rangle\right|} 
\newcommand{\menge}[2]{\big\{{#1}~\big |~{#2}\big\}}

\newcommand{\HH}{\ensuremath{{\mathcal H}}}
\newcommand{\GG}{\ensuremath{{\mathcal G}}}
\newcommand{\supp}{\ensuremath{{\operatorname{supp}}\,}}
\newcommand{\ran}{\ensuremath{{\operatorname{ran}}\,}}
\newcommand{\cran}{\ensuremath{\overline{\operatorname{ran}}\,}}
\newcommand{\spa}{\ensuremath{{\operatorname{span}}\,}}
\newcommand{\spc}{\ensuremath{\overline{\operatorname{span}}\,}}
\newcommand{\emp}{\ensuremath{{\varnothing}}}
\newcommand{\Id}{\ensuremath{\operatorname{Id}}\,}
\newcommand{\RR}{\ensuremath{\mathbb{R}}}
\newcommand{\CC}{\ensuremath{\mathbb{C}}}
\newcommand{\RP}{\ensuremath{\left[0,+\infty\right[}}
\newcommand{\RPP}{\ensuremath{\left]0,+\infty\right[}}
\newcommand{\NN}{\ensuremath{\mathbb N}}
\newcommand{\cart}{\ensuremath{\raisebox{-0.5mm}{\mbox{\LARGE{$\times$}}}}\!}
\newcommand{\ZZ}{\ensuremath{\mathbb Z}}

\newcommand{\exi}{\ensuremath{\exists\,}}
\newcommand{\pinf}{\ensuremath{{+\infty}}}
\newcommand{\minf}{\ensuremath{{-\infty}}}

\newcommand{\zeroun}{\ensuremath{\left]0,1\right[}}   
\newcommand{\Zeroun}{\ensuremath{\left[0,1\right[}}   

\newtheorem{theorem}{Theorem}[section]
\newtheorem{lemma}[theorem]{Lemma}

\newtheorem{corollary}[theorem]{Corollary}
\newtheorem{proposition}[theorem]{Proposition}
\theoremstyle{plain}{\theorembodyfont{\rmfamily}%
}
\theoremstyle{plain}{\theorembodyfont{\rmfamily}%
}
\theoremstyle{plain}{\theorembodyfont{\rmfamily}%
\newtheorem{example}[theorem]{Example}}
\theoremstyle{plain}{\theorembodyfont{\rmfamily}%
\newtheorem{remark}[theorem]{Remark}}
\theoremstyle{plain}{\theorembodyfont{\rmfamily}%
\newtheorem{definition}[theorem]{Definition}}
\theoremstyle{plain}{\theorembodyfont{\rmfamily}%
}

\numberwithin{equation}{section}
\begin{document}
\maketitle

\begin{abstract}
A common problem in applied mathematics is to find a function 
in a Hilbert space with prescribed best approximations from a finite
number of closed vector subspaces. In the present paper we study the 
question of the existence of solutions to such problems. A finite 
family of subspaces is said to satisfy the 
\emph{Inverse Best Approximation Property (IBAP)} if there exists 
a point that admits any selection of points from these subspaces as best 
approximations. We provide various characterizations of the IBAP 
in terms of the geometry of the subspaces. Connections between the IBAP 
and the linear convergence rate of the periodic projection algorithm for 
solving the underlying affine feasibility problem are also established.
The results are applied to problems in harmonic analysis, integral
equations, signal theory, and wavelet frames.
\end{abstract}

\section{Introduction}

A classical problem arising in areas such as harmonic analysis, optics, 
and signal theory is to find a function $x\in L^2(\RR^N)$ with prescribed
values on subsets of the space (or time) and Fourier domains
\cite{Byrn05,Dono89,Havi94,Mele96,Star81,Star87}. In geometrical terms,
this problem can be abstracted into that of finding a function 
possessing prescribed best approximations from two closed vector 
subspaces of 
$L^2(\RR^N)$ \cite{Youl78}. More generally, a broad range of problems 
in applied mathematics can be formulated as follows: given $m$ closed 
vector subspaces $(U_i)_{1\leq i\leq m}$ of a (real or complex)
Hilbert space $\HH$, 
\begin{equation}
\label{e:puertoprincessa-mai2008-1}
\text{find}\;\;x\in\HH\;\;\text{such that}\;\;
(\forall i\in\{1,\ldots,m\})\quad P_ix=u_i,
\end{equation}
where, for every $i\in\{1,\ldots,m\}$, $P_i$ is the (metric)
projector onto $U_i$ and $u_i\in U_i$.
In connection with \eqref{e:puertoprincessa-mai2008-1}, a central 
question is whether a solution exists, irrespective of the choice of 
the prescribed best linear approximations $(u_i)_{1\leq i\leq m}$. 
The main objective of the present paper is to address this question. 

\begin{definition}
\label{d:el-nido2009-03-07}
Let $(U_i)_{1\leq i\leq m}$ be a family of closed vector subspaces of 
$\HH$ and let $(P_i)_{1\leq i\leq m}$ denote their respective
projectors. Then $(U_i)_{1\leq i\leq m}$ satisfies the 
\emph{inverse best approximation property (IBAP)} if 
\begin{equation}
\label{e:manilla2007-9}
\big(\forall (u_i)_{1\leq i\leq m}\in\cart_{i=1}^mU_i\big)
\big(\exi x\in\HH\big)\big(\forall i\in\{1,\ldots,m\})\quad P_ix=u_i.
\end{equation}
Moreover, for every $(u_i)_{1\leq i\leq m}\in\cart_{i=1}^mU_i$, we set
\begin{equation}
\label{e:manilla2008-9}
S(u_1,\ldots,u_m)=\bigcap_{i=1}^m\menge{x\in\HH}{P_ix=u_i},
\end{equation}
and, for every $i\in\{0,\ldots,m-1\}$, 
\begin{equation}
\label{e:palawan-mai2008}
U_{i+}=\sum_{j=i+1}^mU_j\,,\quad P_{i+}=P_{\overline{U_{i+}}}\,,
\quad\text{and}\quad
P_{i+}^\bot=P_{{U_{i+}^\bot}}.
\end{equation}
\end{definition}

The paper is organized as follows. In Section~\ref{sec:2}, we first show 
that the linear independence of the subspaces $(U_i)_{1\leq i\leq m}$ is 
necessary to satisfy the IBAP, but that it is not sufficient in 
infinite dimensional spaces. The main result of Section~\ref{sec:2} is 
Theorem~\ref{t:palawan-mai2008}, which provides various characterizations 
of the IBAP. Several corollaries are derived and, in particular, we 
obtain in Proposition~\ref{p:feasibility} conditions for the 
consistency of affine 
feasibility problems. In Section~\ref{sec:3}, we discuss minimum norm 
solutions and establish connections between the IBAP and the 
rate of convergence of the periodic projection algorithm for solving 
\eqref{e:puertoprincessa-mai2008-1}. Finally, Section~\ref{sec:4} is 
devoted to applications to systems of integral equations, 
constrained moment problems, harmonic analysis, wavelet frames, and 
signal recovery.

\begin{remark}
\label{r:ElNido-mars09}
Since best approximations are well defined for nonempty closed convex
subsets of $\HH$, the IBAP could be considered in this more general
context. However, useful results can be expected to be scarce, even for
two closed convex cones $K_1$ and $K_2$. Indeed, denote the projectors
onto $K_1$ and $K_2$ by $P_1$ and $P_2$, respectively. If $k_1$ is a 
point on the boundary of $K_1$ which is not a support point of $K_1$ (by 
the Bishop-Phelps theorem \cite[Theorem~3.18(i)]{Phel93} support 
points are dense in the boundary of $K_1$), then the only point 
$x\in\HH$ such that $P_1x=k_1$ is $x=k_1$. Therefore, 
there is no point $x\in\HH$ such that $P_1x=k_1$ and $P_2x=k_2$
unless $k_2=P_2k_1$, which means that the IBAP does not hold.
Let us add that, even if every boundary point of $K_1$ is a support 
point (e.g., the interior of $K_1$ is nonempty or $\HH$ is 
finite dimensional), the IBAP can also trivially fail: take for instance 
$\HH=\RR^2$, $K_1=\RP\times\RP$, 
$K_2=\menge{(\beta,-\beta)}{\beta\in\RR}$, 
$k_1=(0,1)$, and $k_2=(1,-1)$.
\end{remark}

Throughout, $\HH$ is a real or complex Hilbert space with scalar product 
$\scal{\cdot}{\cdot}$ and norm $\|\cdot\|$. The distance to a closed
affine subspace $S$ of $\HH$ is denoted by $d_S$, and its projector by 
$P_S$. Moreover, $(U_i)_{1\leq i\leq m}$ is a fixed family of closed 
vector subspaces of $\HH$ with respective projectors 
$(P_i)_{1\leq i\leq m}$. 

\section{Characterizations of the inverse best approximation property}
\label{sec:2}

We first record some useful descriptions of the set of solutions to 
\eqref{e:puertoprincessa-mai2008-1}.

\begin{proposition}
\label{p:carac}
Let $(u_i)_{1\leq i\leq m}\in\cart_{i=1}^mU_i$. Then the following hold.
\begin{enumerate}
\item
\label{p:caraci}
$S(u_1,\ldots,u_m)=\bigcap_{i=1}^m(u_i+U_i^\bot)$.
\item
\label{p:caracii}
Let $x\in S(u_1,\ldots,u_m)$. Then 
$S(u_1,\ldots,u_m)=x+\bigcap_{i=1}^m U_i^{\bot}$.
\end{enumerate}
\end{proposition}
\begin{proof}
\ref{p:caraci}:
Let $x\in\HH$ and $i\in\{1,\ldots,m\}$. The projection theorem asserts 
that $P_ix=u_i$ $\Leftrightarrow$ $x-u_i\in U_i^\bot$
$\Leftrightarrow$ $x\in u_i+U_i^\bot$. Hence, \eqref{e:manilla2008-9} 
yields $x\in S(u_1,\ldots,u_m)$ $\Leftrightarrow$ 
$x\in\bigcap_{i=1}^m(u_i+U_i^\bot)$.

\ref{p:caracii}:
Let $y\in\HH$. By linearity of the operators $(P_i)_{1\leq i\leq m}$,
$y\in S(u_1,\ldots,u_m)$ $\Leftrightarrow$ $(\forall i\in\{1,\ldots,m\})$
$P_i(y-x)=0$ $\Leftrightarrow$ $(\forall i\in\{1,\ldots,m\})$
$y-x\in U_i^\bot$ 
$\Leftrightarrow y\in x+\bigcap_{i=1}^m U_i^{\bot}$.
%
\end{proof}

The main objective of this section is to provide characterizations of 
the inverse best approximation property. Let us start with a necessary 
condition.

\begin{proposition}
\label{p:kimono-ken}
Let $(u_i)_{1\leq i\leq m}\in(\cart_{i=1}^mU_i)
\smallsetminus\{(0,\ldots,0)\}$
be such that $\sum_{i=1}^mu_i=0$. Then $S(u_1,\ldots,u_m)=\emp$. 
\end{proposition}
\begin{proof}
Suppose that $x\in S(u_1,\ldots,u_m)$. Then, for every 
$i\in\{1,\ldots,m\}$,
$u_i=P_ix$ and therefore $\scal{u_i}{x-u_i}=0$, i.e.,
$\|u_i\|^2=\scal{u_i}{x}$. Hence
$0<\sum_{i=1}^m\|u_i\|^2=\scal{\sum_{i=1}^mu_i}{x}=0$, 
and we reach a contradiction.
\end{proof}

\begin{corollary}
\label{c:kimono-ken}
Suppose that $(U_i)_{1\leq i\leq m}$ satisfies the inverse best 
approximation property. Then the subspaces $(U_i)_{1\leq i\leq m}$ are 
linearly independent.
\end{corollary}

As the following example shows, the linear independence of the 
subspaces $(U_i)_{1\leq i\leq m}$ is not sufficient to 
guarantee the inverse best approximation property.

\begin{example}
\label{ex:7}
Suppose that $\HH$ is separable, let $(e_n)_{n\in\NN}$ be an orthonormal 
basis of $\HH$, let $(\alpha_n)_{n\in\NN}$ be a square-summable sequence 
in $\RPP$, and set $(\forall n\in\NN)$ 
$f_n=(e_{2n}+\alpha_ne_{2n+1})/\sqrt{1+\alpha_n^2}$. Set $m=2$,
\begin{equation}
U_1=\spc\{e_{2n}\}_{n\in\NN},\; U_2=\spc\{f_n\}_{n\in\NN},\; u_1=0,\;
\;\text{and}\;\;u_2=\sum_{n\in\NN}\alpha_nf_n.
\end{equation}
Then $U_1\cap U_2=\{0\}$ and $S(u_1,u_2)=\emp$.
\end{example}
\begin{proof}
By construction, $(e_{2n})_{n\in\NN}$ and $(f_n)_{n\in\NN}$ are 
orthonormal bases of $U_1$ and $U_2$, respectively. It follows easily that 
$U_1\cap U_2=\{0\}$. Now suppose that there exists a vector $x\in\HH$ such 
that $_1x=u_1$ and $P_2x=u_2$. Then the identities 
$\sum_{n\in\NN}\scal{x}{e_{2n}}e_{2n}=P_1x=u_1=0$ imply that 
\begin{equation}
\label{e:3-0}
(\forall n\in\NN)\quad\scal{x}{e_{2n}}=0.
\end{equation}
Hence, it results from the identities 
$\sum_{n\in\NN}\alpha_nf_n=u_2=P_2x=\sum_{n\in\NN}\scal{x}{f_n}f_n$ that
\begin{equation}
(\forall n\in\NN)\quad\alpha_n=
\scal{x}{f_n}=\frac{\alpha_n}{\sqrt{1+\alpha_n^2}}\scal{x}{e_{2n+1}}.
\end{equation}
Therefore, $\inf_{n\in\NN}\scal{x}{e_{2n+1}}=\inf_{n\in\NN}
\sqrt{1+\alpha_n^2}=1$, which is impossible.
\end{proof}

The next result states that linear independence is necessary and 
sufficient to obtain an approximate inverse best approximation property. 

\begin{proposition}  
\label{p:denseness}
The following are equivalent.
\begin{enumerate}
\item
\label{p:densenessi}
The subspaces $(U_i)_{1\leq i\leq m}$ are linearly independent.
\item
\label{p:densenessii}
For every $(u_i)_{1\leq i\leq m}\in\cart_{i=1}^mU_i$ and every
$\varepsilon\in\RPP$, there exists $x\in\HH$ such that
\begin{equation}
\label{e:23octobre2008}
\max_{1\leq i\leq m}\|P_ix-u_i\|\leq\varepsilon.
\end{equation}
\end{enumerate}
\end{proposition}
\begin{proof}
Set $V=\menge{(P_ix)_{1\leq i\leq m}}{x\in\HH}$ and let $W$ be the
orthogonal complement of $V$ in the Hilbert direct sum
$\bigoplus_{i=1}^mU_i$.

\ref{p:densenessi}$\Rightarrow$\ref{p:densenessii}:
Take $(u_i)_{1\leq i\leq m}\in W$ and set $x=\sum_{i=1}^{m}u_i$. Then 
$\sum_{i=1}^{m}\scal{u_i}{x}=\sum_{i=1}^{m}\scal{u_i}{P_ix}=0$, which
implies that $\|x\|^2=\sum_{i=1}^{m}\scal{u_i}{x}=0$. Hence $x=0$ and, 
in view of the assumption of independence, we conclude that 
$(\forall i\in\{1,\ldots,m\})$ $u_i=0$.
Therefore, $V$ is dense in $\bigoplus_{i=1}^mU_i$.

\ref{p:densenessii}$\Rightarrow$\ref{p:densenessi}:
Take $(u_i)_{1\leq i\leq m}\in\cart_{i=1}^mU_i$ such that
$\sum_{i=1}^{m}u_i=0$, take $\varepsilon\in\RPP$, and take 
$x\in\HH$ such that \eqref{e:23octobre2008} holds. Then
$\sum_{i=1}^m\scal{u_i}{P_ix}=\sum_{i=1}^m\scal{u_i}{x}=0$ 
and therefore
\begin{align}
\sum_{i=1}^m\|u_i\|^2
&=\sum_{i=1}^m\|u_i-P_ix\|^2+2 \text{Re} \sum_{i=1}^m\scal{u_i-P_ix}{P_ix}
+\sum_{i=1}^m\|P_ix\|^2\nonumber\\
&=\sum_{i=1}^m\|u_i-P_ix\|^2-\sum_{i=1}^m\|P_ix\|^2\nonumber\\
&\leq m\varepsilon^2.
\end{align}
Hence, $(\forall i\in\{1,\ldots,m\})$ $u_i=0$. 
\end{proof}

In order to provide characterizations of the inverse best 
approximation property, we require the following tools.

\begin{definition}{\rm\cite[Definition~9.4]{Deut01}}
Let $U$ and $V$ be closed vector subspaces of $\HH$.
The angle determined by $U$ and $V$ is the real 
number in $[0,\pi/2]$ the cosine of which is given by
\begin{equation}
\label{e:angol}
{\mathsf c}(U,V)
=\sup\menge{\abscal{x}{y}}{x\in U\cap(U\cap V)^{\bot},\:
y\in V\cap(U\cap V)^{\bot},\:\|x\|\leq 1,\:\|y\|\leq 1}.
\end{equation}
\end{definition}

\begin{lemma}
\label{l:palawan-mai2008}
Let $U$ and $V$ be closed vector subspaces of $\HH$, let 
$u\in U$, let $v\in V$, and set $S=(u+U^\bot)\cap(v+V^\bot)$.
Then the following hold.
\begin{enumerate}
\item
\label{l:palawan-mai2008i}
Let $x\in S$. Then $S=P_{\overline{U+V}}\,x+(U^\bot\cap V^\bot)$.
\item
\label{l:palawan-mai2008ii}
Suppose that $\|P_UP_V\|<1$ and set
\begin{equation}
\label{e:teryaki-boy3}
z={\overline{u}}+{\overline{v}},\quad\text{where}\quad
\begin{cases}
{\overline{u}}=(\Id-P_UP_V)^{-1}(u-P_Uv)\\
{\overline{v}}=(\Id-P_VP_U)^{-1}(v-P_Vu).
\end{cases}
\end{equation}
Then the following hold.
\begin{enumerate}
\item
\label{l:palawan-mai2008iia}
$S\neq\emp$. 
\item
\label{l:palawan-mai2008iib}
$z=P_S\,0$.
\end{enumerate}
\end{enumerate}
\end{lemma}
\begin{proof}
\ref{l:palawan-mai2008i}: 
As in Proposition~\ref{p:carac}, we can write 
$S=x+(U^{\bot}\cap V^{\bot})$. Hence, since 
$(\overline{U+V})^\bot=(U+V)^\bot=U^\bot\cap V^\bot$, we get
$S=x+(U^{\bot}\cap V^{\bot})
=P_{(U^{\bot}\cap V^{\bot})^\bot}x+(U^{\bot}\cap V^{\bot})
=P_{\overline{U+V}}\,x+(U^{\bot}\cap V^{\bot})$.

\ref{l:palawan-mai2008ii}:
These properties are known (see for instance \cite[Item~3.B) p.~91]{Havi94}
and \cite[Section~5 on pp.~92--93]{Havi94}, respectively); 
we provide short alternative proofs for completeness.

\ref{l:palawan-mai2008iia}:
Let $u\in U$ and $v\in V$. Since $P_U$ and $P_V$ are self-adjoint,
$\|P_VP_U\|=\|(P_VP_U)^*\|=\|P_U^*P_V^*\|=\|P_UP_V\|<1$, and the vectors
${\overline{u}}$ and ${\overline{v}}$ are therefore well defined. 
Moreover, it follows from the identity
${\overline{u}}=\sum_{j\in\NN}(P_UP_V)^j(u-P_Uv)$ that 
${\overline{u}}\in U$ and therefore that
$P_U{\overline{u}}=\overline{u}$. On the other hand, the second equality 
in the right-hand side of \eqref{e:teryaki-boy3} yields
\begin{align}
\label{e:teryaki-boy1}
P_U{\overline{v}}
&=P_U\bigg(\sum_{j\in\NN}(P_V P_U)^j(v-P_Vu)\bigg)\nonumber\\
&=(\Id-P_UP_V)^{-1}(P_Uv-P_UP_Vu)\nonumber\\
&=(\Id-P_UP_V)^{-1}\big((\Id-P_UP_V)u-(u-P_Uv)\big)\nonumber\\
&=u-\overline{u}.
\end{align}
Thus, $P_Uz=P_U({\overline{u}}+{\overline{v}})=
{\overline{u}}+P_U{\overline{v}}=u$.
Likewise, $P_V{\overline{v}}=\overline{v}$ and
$P_V\overline{u}=v-\overline{v}$, which implies that 
$P_Vz=P_V({\overline{u}}+{\overline{v}})=
P_V{\overline{u}}+{\overline{v}}=v$. Altogether, $z\in S$.

\ref{l:palawan-mai2008iib}:
As seen above, $z\in S$, ${\overline{u}}\in U$, and 
${\overline{v}}\in V$. Now let $x\in S$.
As in Proposition~\ref{p:carac}\ref{p:caracii}, we can write
$x=z+w=\overline{u}+\overline{v}+w$, for some 
$w\in U^\bot\cap V^\bot$. Hence,
$\|x\|^2=\|z\|^2+2\text{Re}\scal{\overline{u}}{w}+2\text{Re}
\scal{\overline{v}}{w}+\|w\|^2=\|z\|^2+\|w\|^2\geq\|z\|^2$.
\end{proof}

We can now provide various characterizations of the inverse best 
approximation property (the notation \eqref{e:palawan-mai2008} will 
be used repeatedly). 

\begin{theorem}
\label{t:palawan-mai2008}
The following are equivalent.
\begin{enumerate}
\item
\label{t:palawan-mai2008i}
$(U_i)_{1\leq i\leq m}$ satisfies the inverse best approximation property.
\item
\label{t:palawan-mai2008ii}
$(\forall i\in\{1,\ldots,m-1\})(\forall u_i\in U_i)(\exi x\in\HH)$ 
$u_i=P_ix$ and $(\forall j\in\{i+1,\ldots,m\})$ $P_jx=0$.
\item
\label{t:palawan-mai2008iii}
$(\forall i\in\{1,\ldots,m-1\})$ $P_i(U_{i+}^\bot)=U_i$.
\item
\label{t:palawan-mai2008iv}
$(\forall i\in\{1,\ldots,m-1\})$ $U_i^\bot+U_{i+}^\bot=\HH$.
\item
\label{t:palawan-mai2008x}
The subspaces $(U_i)_{1\leq i\leq m}$ are linearly independent and
$(\forall i\in\{1,\ldots,m-1\})(\exi\gamma_i\in\RPP)$ 
$d_{U_i^\bot\cap U^\bot_{i+}}\leq\gamma_i\big(
d_{U^\bot_i}+d_{U^\bot_{i+}}\big)$.
\item
\label{t:palawan-mai2008v}
The subspaces $(U_i)_{1\leq i\leq m}$ are linearly independent and
$(\forall i\in\{1,\ldots,m-1\})$ $U_i+U_{i+}$ is closed.
\item
\label{t:palawan-mai2008vi}
The subspaces $(U_i)_{1\leq i\leq m}$ are linearly independent and, for
every $i\in\{1,\ldots,m-1\}$, ${\mathsf c}(U_i,U_{i+})<1$.
\item
\label{t:palawan-mai2008vii}
$(\forall i\in\{1,\ldots,m-1\})(\exi\gamma_i\in\left[1,\pinf\right[)
(\forall u_i\in U_i)$ $\|u_i\|\leq\gamma_i\|P_{i+}^\bot u_i\|$.
\item
\label{t:palawan-mai2008vii'}
$(\forall i\in\{1,\ldots,m-1\})(\exi\gamma_i\in\left[2,\pinf\right[)
(\forall x\in\HH)$ $\|x\|\leq\gamma_i(\|P_i^\bot x\|
+\|P_{i+}^\bot x\|)$.
\item
\label{t:palawan-mai2008viii}
$(\forall i\in\{1,\ldots,m-1\})$ $\|P_iP_{i+}\|<1$.
\end{enumerate}
\end{theorem}
\begin{proof}
\ref{t:palawan-mai2008i}$\Rightarrow$\ref{t:palawan-mai2008ii}: Clear.

\ref{t:palawan-mai2008ii}$\Rightarrow$\ref{t:palawan-mai2008iii}: 
Let $i\in\{1,\ldots,m-1\}$. It is clear that 
$P_i(U_{i+}^\bot)\subset U_i$. Conversely, let $u_i\in U_i$. By 
assumption, there exists $x\in\bigcap_{j=i+1}^mU_j^\bot=U_{i+}^\bot$ 
such that $u_i=P_ix$. In other words, $U_i\subset P_i(U_{i+}^\bot)$. 
Altogether, $P_i(U_{i+}^\bot)=U_i$. 

\ref{t:palawan-mai2008iii}$\Rightarrow$\ref{t:palawan-mai2008iv}:
Let $i\in\{1,\ldots,m-1\}$.
We have 
\begin{equation}
\HH=U_i^\bot+U_i=U_i^\bot+P_i(U_{i+}^\bot)
=U_i^\bot+\bigcup_{v\in U_{i+}^\bot}(v-P_{U_i^\bot}v)
=U_i^\bot+\bigcup_{v\in U_{i+}^\bot}v=U_i^\bot+U_{i+}^\bot.
\end{equation}

\ref{t:palawan-mai2008iv}$\Rightarrow$\ref{t:palawan-mai2008x}: 
Let $i\in\{1,\ldots,m-1\}$. We have
\begin{equation}
U_i\cap U_{i+}=(U_i^\bot+U_{i+}^\bot)^\bot=\HH^\bot=\{0\}.
\end{equation}
This shows  the independence claim. Moreover, since 
$U_i^\bot+U_{i+}^\bot=\HH$ is closed, the inequality on the distance
functions follows from \cite[Corollaire~II.9]{Brez93}.

\ref{t:palawan-mai2008x}$\Rightarrow$\ref{t:palawan-mai2008v}: 
Let $i\in\{1,\ldots,m-1\}$. It follows from 
\cite[Remarque~7~p.~22]{Brez93} (see also
\cite[Proposition~5.16]{Baus96}) that $U_i^\bot+U_{i+}^\bot$ is closed.
In turn, since \cite[Th\'eor\`eme~II.15]{Brez93} asserts that
$U_i^{\bot\bot}+U^{\bot\bot}_{i+}$ is closed, we deduce that
\begin{equation}
\label{e:brezis82}
U_i+\overline{U_{i+}}~\text{is closed}.
\end{equation}
It remains to show that $U_{i+}$ is closed.
If $i=m-1$, $U_{i+}=U_m$ is closed. On the other hand, if 
$i\in\{2,\ldots,m-1\}$ and $U_{i+}$ is closed, we deduce from
\eqref{e:brezis82} that
$U_{(i-1)+}=U_i+U_{i+}=U_i+\overline{U_{i+}}$ is closed.

\ref{t:palawan-mai2008v}$\Rightarrow$\ref{t:palawan-mai2008vi}: 
Let $i\in\{1,\ldots,m-1\}$. Then
$U_{i+}$ and $U_i+U_{i+}$ are closed and it follows from 
\cite[Theorem~9.35]{Deut01} that ${\mathsf c}(U_i,U_{i+})<1$.

\ref{t:palawan-mai2008vi}$\Rightarrow$\ref{t:palawan-mai2008vii}:
Let $i\in\{1,\ldots,m-1\}$ and let $u_i\in U_i$. Then
\eqref{e:angol} yields
\begin{align}
\label{e:3nov2008-1}
\|u_i\|^2
&=\|P_{i+}^\bot u_i\|^2+\|P_{i+}u_i\|^2\nonumber\\
&=\|P_{i+}^\bot u_i\|^2+\scal{u_i}{P_{i+}u_i}\nonumber\\
&\leq\|P_{i+}^\bot u_i\|^2+{\mathsf c}(U_i,U_{i+})\|u_i\|\,
\|P_{i+}u_i\|\nonumber\\
&\leq\|P_{i+}^\bot u_i\|^2+{\mathsf c}(U_i,U_{i+})\|u_i\|^2. 
\end{align}
Hence,
$\|P^\bot_{i+}u_i\|^2\geq(1-{\mathsf c}(U_i,U_{i+}))\|u_i\|^2$.

\ref{t:palawan-mai2008vii}$\Rightarrow$\ref{t:palawan-mai2008vii'}:
Let $i\in\{1,\ldots,m-1\}$ and let $x\in\HH$. There exists
$\gamma\in\left[1,\pinf\right[$ such that
\begin{align}
\label{e:3nov2008-2}
\|x\|
&\leq\|P_ix\|+\|P^\bot_ix\|\nonumber\\
&\leq\gamma\|P_{i+}^\bot P_ix\|+\|P^\bot_ix\|\nonumber\\
&\leq\gamma(\|P_{i+}^\bot x\|+\|P_{i+}^\bot P^\bot_ix\|)
+\|P^\bot_ix\|\nonumber\\
&\leq\gamma\|P_{i+}^\bot x\|+(1+\gamma)\|P^\bot_ix\|.
\end{align}

\ref{t:palawan-mai2008vii'}$\Rightarrow$\ref{t:palawan-mai2008viii}:
Let $i\in\{1,\ldots,m-1\}$ and let $x\in\HH$. There exists
$\gamma\in\left[2,\pinf\right[$ such that
\begin{align}
\|P_ix\|^2
&=\|P_{i+}P_ix\|^2+\|P_{i+}^\bot P_ix\|^2\nonumber\\
&=\|P_{i+}P_ix\|^2+(\|P_{i+}^\bot P_ix\|+\|P_i^\bot P_ix\|)^2\nonumber\\
&\geq\|P_{i+}P_ix\|^2+\gamma^{-2}\|P_ix\|^2. 
\end{align}
Therefore
$\|P_{i+}P_ix\|^2\leq(1-\gamma^{-2})\|P_ix\|^2\leq(1-\gamma^{-2})\|x\|^2$.
Hence $\|P_{i+}P_i\|<1$ and, in turn, 
$\|P_iP_{i+}\|=\|P_i^*P_{i+}^*\|=\|(P_{i+}P_i)^*\|=\|P_{i+}P_i\|<1$.

\ref{t:palawan-mai2008viii}$\Rightarrow$\ref{t:palawan-mai2008i}:
Fix $(u_i)_{1\leq i\leq m}\in\cart_{i=1}^mU_i$ and set
$(\forall i\in\{0,\ldots,m-1\}$ 
$S_i=\bigcap_{j=i+1}^m(u_i+U_i^\bot)$.
Let us show by induction that
\begin{equation}
\label{e:mary'scottage-mai2008}
(\forall i\in\{0,\ldots,m-2\})\quad 
S_i\neq\emp\quad\text{and}\quad
(\forall x_i\in S_i)\quad S_i=P_{i+}x_i+U_{i+}^\bot.
\end{equation}
First, let us set $i=m-2$. Since, by assumption $\|P_{m-1}P_m\|<1$, it
follows from Lemma~\ref{l:palawan-mai2008}\ref{l:palawan-mai2008iia} that 
$S_{m-2}\neq\emp$. Moreover, we 
deduce from Lemma~\ref{l:palawan-mai2008}\ref{l:palawan-mai2008i} that,
for every $x_{m-2}\in S_{m-2}$,
\begin{equation}
S_{m-2}=P_{\overline{U_{m-1}+U_m}}\,x_{m-2}+(U_{m-1}^\bot\cap
U_m^\bot)=P_{(m-2)+}x_{m-2}+U_{(m-2)+}^\bot.
\end{equation}
Next, suppose that \eqref{e:mary'scottage-mai2008} is true for some 
$i\in\{1,\ldots,m-2\}$ and let $x_i\in S_i$.
Then, using Lemma~\ref{l:palawan-mai2008}\ref{l:palawan-mai2008i},
we obtain 
\begin{equation}
\label{e:rains} 
S_{i-1}
=(u_i+U_i^\bot)\cap S_i
=(u_i+U_i^\bot)\cap\big(P_{i+}x_i+U_{i+}^\bot).
\end{equation}
Since, by assumption $\|P_iP_{i+}\|<1$, it follows from 
Lemma~\ref{l:palawan-mai2008}\ref{l:palawan-mai2008iia} that
$S_{i-1}\neq\emp$. Now, let $ x_{i-1}\in S_{i-1}$.
Combining  \eqref{e:rains} and Lemma~\ref{l:palawan-mai2008}
\ref{l:palawan-mai2008i}, we obtain
\begin{equation}
S_{i-1}=P_{\overline{U_i+U_{i+}}}x_{i-1} +
\left( U_i^{\bot}\cap U_{i+}^{\bot}\right)
= P_{(i-1)+} x_{i-1} + U_{(i-1)+}^{\bot}.
\end{equation}
This proves by
induction that \eqref{e:mary'scottage-mai2008} is true. For 
$i=0$, we thus obtain $S_0=\bigcap_{j=1}^m(u_j+U_j^\bot)\neq\emp$.
In view of Proposition~\ref{p:carac}\ref{p:caraci}, the proof is 
complete.
\end{proof}

An immediate application of Theorem~\ref{t:palawan-mai2008} concerns 
the area of affine feasibility problems
\cite{Kruk06,Butn01,Byrn05,Proc93,Kosm91,Star87}. Given a family of
closed affine subspaces $(S_i)_{1\leq i\leq m}$ of $\HH$, the problem
is to 
\begin{equation}
\label{e:palawan-mars2009-1}
\text{find}\;\;x\in\bigcap_{i=1}^mS_i.
\end{equation}
In applications, a key issue is whether this problem is consistent in 
the sense that it admits a solution. Our next proposition gives a 
sufficient condition for consistency. First, we recall a standard fact.

\begin{lemma}
\label{l:45}
Let $S$ be a closed affine subspace of $\HH$, let $V=S-S$ be the closed
vector subspace parallel to $S$, and let $y\in S$. Then $S=y+V$ and
$(\forall x\in\HH)$ $P_Sx=y+P_V(x-y)$.
\end{lemma}

\begin{proposition}
\label{p:feasibility}
Let $(S_i)_{1\leq i\leq m}$ be closed affine subspaces of $\HH$ 
and suppose that $(U_i)_{1\leq i\leq m}$ are the orthogonal 
complements of their respective parallel vector subspaces. If 
$(U_i)_{1\leq i\leq m}$ satisfies the inverse best approximation 
property (in particular, if any of 
properties~{\rm\ref{t:palawan-mai2008ii}--\ref{t:palawan-mai2008viii}} 
in Theorem~{\rm\ref{t:palawan-mai2008}} holds), then the affine 
feasibility problem \eqref{e:palawan-mars2009-1} is consistent.
\end{proposition}
\begin{proof}
For every $i\in\{1,\ldots,m\}$, let $a_i\in S_i$, and set 
$V_i=S_i-S_i$ and $u_i=P_ia_i$. Then, by Lemma~\ref{l:45},
$(\forall i\in\{1,\ldots,m\})$
$S_i=a_i+V_i=a_i+U_i^{\bot}=u_i+U_i^\bot$. Thus, 
\begin{equation}
\bigcap_{i=1}^mS_i=\bigcap_{i=1}^m(u_i+U_i^{\bot}),
\end{equation}
and it follows from Proposition~\ref{p:carac}\ref{p:caraci} that
\eqref{e:palawan-mars2009-1} is consistent if
$(U_i)_{1\leq i\leq m}$ satisfies the IBAP. 
\end{proof}

\begin{remark}
The converse to Proposition~\ref{p:feasibility} fails. For instance, 
let $S_1$ and $S_2$ be distinct intersecting lines in $\HH=\RR^3$.
Then $U_1=(S_1-S_1)^\bot$ and $U_2=(S_2-S_2)^\bot$ are two-dimensional
planes and they are therefore linearly dependent. Hence, the IBAP cannot
hold by virtue of Corollary~\ref{c:kimono-ken}.
\end{remark}

In the case of two subspaces, Theorem~\ref{t:palawan-mai2008} yields 
simpler conditions.
\begin{corollary}
\label{c:palawan-mai2008}
The following are equivalent.
\begin{enumerate}
\item
\label{c:palawan-mai2008i}
$(U_1,U_2)$ satisfies the inverse best approximation property.
\item
\label{c:palawan-mai2008ii}
$(\forall u_1\in U_1)$ $S(u_1,0)\neq\emp$.
\item
\label{c:palawan-mai2008iii}
$P_1(U_2^\bot)=U_1$.
\item
\label{c:palawan-mai2008iv}
$U_1^\bot+U_2^\bot=\HH$.
\item
\label{c:palawan-mai2008x}
$U_1\cap U_2=\{0\}$ and $(\exi\gamma\in\RPP)$ 
$d_{U^\bot_1\cap U^\bot_2}\leq\gamma\big(d_{U^\bot_1}+d_{U^\bot_2}\big)$.
\item
\label{c:palawan-mai2008v}
$U_1\cap U_2=\{0\}$ and $U_1+U_2$ is closed.
\item
\label{c:palawan-mai2008vi}
$U_1\cap U_2=\{0\}$ and ${\mathsf c}(U_1,U_2)<1$.
\item
\label{c:palawan-mai2008vii}
$(\exi\gamma\in\left[1,\pinf\right[)
(\forall u_1\in U_1)$ $\|u_1\|\leq\gamma\|P_2^\bot u_1\|$.
\item
\label{c:palawan-mai2008vii+}
$(\exi\gamma\in\left[2,\pinf\right[)
(\forall x\in\HH)$ $\|x\|\leq\gamma(\|P_1^\bot x\|+\|P_2^\bot x\|)$.
\item
\label{c:palawan-mai2008viii}
$\|P_1P_2\|<1$.
\end{enumerate}
\end{corollary}

\begin{remark}
Corollary~\ref{c:palawan-mai2008} provides necessary and 
sufficient conditions for the existence of solutions to 
\eqref{e:puertoprincessa-mai2008-1} when $m=2$. The implication 
\ref{c:palawan-mai2008vii+}$\Rightarrow$\ref{c:palawan-mai2008i} 
appears in \cite[Item~3.B)~p.~91]{Havi94}, the equivalences 
\ref{c:palawan-mai2008v}$\Leftrightarrow$\ref{c:palawan-mai2008vii}%
$\Leftrightarrow$\ref{c:palawan-mai2008vii+}%
$\Leftrightarrow$\ref{c:palawan-mai2008viii} appear in 
\cite[Item~1.A)~p.~88]{Havi94}, and the equivalences 
\ref{c:palawan-mai2008iii}$\Leftrightarrow$\ref{c:palawan-mai2008iv}%
$\Leftrightarrow$\ref{c:palawan-mai2008viii}
appear in \cite[Lemma on p.~201]{Niko86}. 
\end{remark}

As consequences of Theorem~\ref{t:palawan-mai2008}, we can now 
describe scenarii in which the necessary condition established in 
Corollary~\ref{c:kimono-ken} is also sufficient.

\begin{corollary}
\label{c:paris-octobre2008}
Suppose that the closed vector subspaces $(U_i)_{1\leq i\leq m}$ are 
linearly independent, that $\|P_{m-1}P_m\|<1$ and that, for every
$i\in\{1,\ldots,m-2\}$, $U_i$ is finite dimensional or finite 
codimensional. Then $(U_i)_{1\leq i\leq m}$ satisfies the inverse 
best approximation property.
\end{corollary}
\begin{proof}
In view of the equivalence 
\ref{t:palawan-mai2008i}$\Leftrightarrow$\ref{t:palawan-mai2008vi}
in Theorem~\ref{t:palawan-mai2008}, it is enough to show that
$(\forall i\in\{1,\ldots,m-1\})$ ${\mathsf c}(U_i,U_{i+})<1$.
For $i=m-1$, since $\|P_iP_{i+}\|=\|P_{m-1}P_m\|<1$, we derive from 
the implication 
\ref{c:palawan-mai2008viii}$\Rightarrow$\ref{c:palawan-mai2008vi} in
Corollary~\ref{c:palawan-mai2008} that ${\mathsf c}(U_i,U_{i+})<1$.
Now suppose that, for some $i\in\{2,\ldots,m-1\}$, 
${\mathsf c}(U_i,U_{i+})<1$. Using to the implication
\ref{c:palawan-mai2008vi}$\Rightarrow$\ref{c:palawan-mai2008v} in
Corollary~\ref{c:palawan-mai2008}, we deduce that
$U_{(i-1)+}=U_i+U_{i+}$ is closed. In turn, since 
$U_{i-1}$ is finite or cofinite dimensional, it follows from 
\cite[Corollary~9.37]{Deut01} that ${\mathsf c}(U_{(i-1)},U_{(i-1)+})<1$, 
which completes the proof by induction.
\end{proof}

\begin{corollary}
\label{c:puertoprincessa-mai2008}
Suppose that the closed vector subspaces $(U_i)_{1\leq i\leq m}$ are 
linearly independent and that, for every $i\in\{1,\ldots,m-1\}$, $U_i$ 
is finite dimensional or finite codimensional. Then 
$(U_i)_{1\leq i\leq m}$ satisfies the inverse best approximation property.
\end{corollary}
\begin{proof}
Since $U_{m-1}$ is finite dimensional or finite codimensional, it follows 
from \cite[Corollary~9.37]{Deut01} and the implication
\ref{c:palawan-mai2008vi}$\Rightarrow$\ref{c:palawan-mai2008viii} in
Corollary~\ref{c:palawan-mai2008} that $\|P_{m-1}P_m\|<1$. Hence, 
the claim follows from Corollary~\ref{c:paris-octobre2008}.
\end{proof}

\begin{example}
\label{prob:72}
Let $V$ be a closed vector subspace of $\HH$ and let 
$(v_i)_{1\leq i\leq m-1}$ be linearly independent vectors
such that $V^\bot\cap \spa\{v_i\}_{1\leq i\leq m-1}=\{0\}$. 
Then, for every $(\eta_i)_{1\leq i\leq m-1}\in\CC^{m-1}$, 
the constrained moment problem 
\begin{equation}
\label{e:p72}
x\in V\quad\text{and}\quad
(\forall i\in\{1,\ldots,m-1\})\quad\scal{x}{v_i}=\eta_i
\end{equation}
admits a solution.
\end{example}
\begin{proof}
This is a special case of Corollary~\ref{c:puertoprincessa-mai2008}, 
where $U_m=V^\bot$, $u_m=0$, and, for every $i\in\{1,\ldots,m-1\}$,
$U_i=\spa\{v_i\}$ and $u_i=\eta_iv_i/\|v_i\|^2$. 
\end{proof}

\begin{corollary}
\label{c:puertoprincessa-mai2008-2}
Suppose that the subspaces $(U_i)_{1\leq i\leq m}$ are linearly
independent and that $\HH$ is finite dimensional. Then 
$(U_i)_{1\leq i\leq m}$ satisfies the inverse best approximation property.
\end{corollary}

The above results pertain to the existence of solutions to
\eqref{e:puertoprincessa-mai2008-1}. We conclude this section
with a uniqueness result that follows at once from 
Proposition~\ref{p:carac}\ref{p:caracii}.

\begin{proposition}
\label{p:23}
Let $(u_i)_{1\leq i\leq m}\in\cart_{i=1}^mU_i$. Then
\eqref{e:puertoprincessa-mai2008-1} has at most one solution 
if and only if $\bigcap_{i=1}^mU_i^\bot=\{0\}$.
\end{proposition}

Combining Theorem~\ref{t:palawan-mai2008} and Proposition~\ref{p:23} 
yields conditions for the existence of unique solutions to
\eqref{e:puertoprincessa-mai2008-1}. Here is an example in which $m=2$.

\begin{example}
\label{ex:unique1}
The following are equivalent.
\begin{enumerate}
\item
\label{ex:unique1i}
For every $u_1\in U_1$ and $u_2\in U_2$, $S(u_1,u_2)$ is a singleton.
\item
\label{ex:unique1ii}
$U_1^\bot+U_2^\bot=\HH$ and $U_1^{\bot}\cap U_2^{\bot}=\{0\}$. 
\end{enumerate}
\end{example}
\begin{proof}
Existence follows from the implication 
\ref{c:palawan-mai2008iv}$\Rightarrow$\ref{c:palawan-mai2008i} in 
Corollary~\ref{c:palawan-mai2008}, and uniqueness from 
Proposition~\ref{p:23}.
\end{proof}

\section{IBAP and the periodic projection algorithm}
\label{sec:3}

If $(U_i)_{1\leq i\leq m}$ satisfies the IBAP, then
\eqref{e:puertoprincessa-mai2008-1} will in general admit infinitely
many solutions (see Proposition~\ref{p:23}) and it is of interest to 
identify specific solutions such as those of minimum norm.

\begin{proposition}
\label{p:24}
Suppose that $(U_i)_{1\leq i\leq m}$ satisfies the inverse best
approximation property, let 
$(u_i)_{1\leq i\leq m}\in\cart_{i=1}^mU_i$, and, 
for every $i\in\{1,\ldots,m-1\}$, set
\begin{eqnarray}
\label{e:snakeisland-mai2008}
T_i\colon U_{i+}&\to& U_i+U_{i+}\nonumber\\
v~&\mapsto &(\Id-P_iP_{i+})^{-1}(u_i-P_iv)+
(\Id-P_{i+}P_i)^{-1}(v-P_{i+}u_i).
\end{eqnarray}
Define recursively $\overline{x}_m=u_m$ and 
$(\forall i\in\{m-1,\ldots,1\})$ $\overline{x}_i=T_i\overline{x}_{i+1}$.
Then, for every $i\in\{1,\ldots,m\}$,
\begin{equation}
\label{e:minimal}
\overline{x}_i=P_{S_i}0,\quad\text{where}\quad
S_i=\bigcap_{j=i}^m(u_j+U_j^\bot).
\end{equation}
In particular, $\overline{x}_1=P_{S(u_1,\ldots,u_m)}0$ is the minimal
norm solution to \eqref{e:puertoprincessa-mai2008-1}.
\end{proposition}
\begin{proof} 
Let $i\in\{1,\ldots,m-1\}$. We first observe that the operator 
$T_i$ is well defined since the implication
\ref{t:palawan-mai2008i}$\Rightarrow$\ref{t:palawan-mai2008viii}
in Theorem~\ref{t:palawan-mai2008} yields
$\|P_iP_{i+}\|=\|P_{i+}P_i\|<1$. Moreover,
the expansions $(\Id-P_iP_{i+})^{-1}=\sum_{j\in\NN}(P_iP_{i+})^j$
and $(\Id-P_{i+}P_i)^{-1}=\sum_{j\in\NN}(P_{i+}P_i)^j$ imply that
its range is indeed contained in $U_i+U_{i+}$. Thus, 
$\overline{x}_i$ is a well defined point in
$U_i+U_{i+}=U_{(i-1)+}$.

To prove \eqref{e:minimal}, we proceed by induction. First, for $i=m$, 
since $u_m\in U_m$, we obtain at once 
\begin{equation}
\overline{x}_{i}=u_m=P_{(u_m+U_m^\bot)}0=P_{S_i}0.
\end{equation}
Now, suppose that \eqref{e:minimal} is true for some 
$i\in\{2,\ldots,m\}$. By definition,
\begin{equation}
\overline{x}_{i-1}=
(\Id-P_{i-1}P_{(i-1)+})^{-1}(u_{i -1}-P_{i-1}\overline{x}_i)+
(\Id-P_{(i-1)+}P_{i-1})^{-1}(\overline{x}_i-P_{(i-1)+}u_{i-1}).
\end{equation}
Since $\overline{x}_i\in U_{(i-1)+}$ and $u_{i-1}\in U_{i-1}$,
Lemma~\ref{l:palawan-mai2008}\ref{l:palawan-mai2008iib} asserts that 
$\overline{x}_{i-1}$ is the element of minimal norm in
$(u_{i-1}+U^{\bot}_{i-1})\cap(\overline{x}_i+U^{\bot}_{(i-1)+})$. 
On the other hand since, by \eqref{e:minimal}, 
$\overline{x}_i\in\bigcap_{j=i}^m(u_j+U_j^\bot)$, we derive from
\eqref{e:palawan-mai2008} that, as in Proposition~\ref{p:carac},
\begin{equation}
\overline{x}_i+U^{\bot}_{(i-1)+}
=\overline{x}_i+\Bigg(\sum_{j=i}^mU_j\Bigg)^\bot
=\overline{x}_i+\bigcap_{j=i}^mU_j^\bot
=\bigcap_{j=i}^m(u_j+U_j^\bot).
\end{equation}
As a result, $\overline{x}_{i-1}$ is the element of minimum norm in
\begin{equation}
(u_{i-1}+U^{\bot}_{i-1})\cap\bigcap_{j=i}^m(u_j+U_j^\bot).
\end{equation}
In other words, $\overline{x}_{i-1}=P_{S_{i-1}}0$, which
completes the proof.
\end{proof}

Conceptually, Proposition~\ref{p:24} provides a finite 
recursion for computing the minimal norm solution $\overline{x}_1$ to 
\eqref{e:puertoprincessa-mai2008-1} for a given selection of
vectors $(u_i)_{1\leq i\leq m}\in\cart_{i=1}^mU_i$. This scheme
is in general not of direct numerical use since it requires the 
inversion of operators in \eqref{e:snakeisland-mai2008}. However,
minimal norm solutions and, more generally, best approximations from
the solution set of \eqref{e:puertoprincessa-mai2008-1} can be 
computed iteratively via projection methods. Indeed, for every 
$r\in\HH$ and $(u_i)_{1\leq i\leq m}\in\cart_{i=1}^mU_i$, let us 
denote by $B(r;u_1,\ldots,u_m)$ the best approximation to $r$ 
from $S(u_1,\ldots,u_m)$, i.e., by 
Proposition~\ref{p:carac}\ref{p:caraci}, 
\begin{equation}
\label{e:diliman2009-03-12a}
B(r;u_1,\ldots,u_m)=P_{S(u_1,\ldots,u_m)}r=
P_{\bigcap_{i=1}^m(u_i+U_i^\bot)}r.
\end{equation}
A standard numerical method for computing $B(r;u_1,\ldots,u_m)$ is the 
periodic projection algorithm 
\begin{equation}
\label{e:el-nido2009-03-10a}
x_0=r\quad\text{and}\quad(\forall n\in\NN)\quad
x_{n+1}=Q_1\cdots Q_mx_n
\end{equation}
where, for every $i\in\{1,\ldots,m\}$, $Q_i$ is the projector onto
$u_i+U_i^\bot$, i.e., 
\begin{equation}
\label{e:puertoprincessa-mars2009-1}
Q_i=P_{u_i+U_i^\bot}\colon x\mapsto u_i+x-P_ix.
\end{equation}
This algorithm is rooted in the classical work of Kaczmarz
\cite{Kacz37} and von Neumann \cite{Vonn49}. Although it has been 
generalized in various directions \cite{Baus96,Kruk06,Butn01,Jamo97}, 
it is still widely used due to its simplicity and ease of implementation. 
If $S(u_1,\ldots,u_m)\neq\emp$, the sequence $(x_n)_{n\in\NN}$ 
generated by \eqref{e:el-nido2009-03-10a} converges strongly to 
$B(r;u_1,\ldots,u_m)$. 
If $u_i\equiv 0$, this result was first established by von Neumann 
\cite{Vonn49} for $m=2$ and extended by Halperin \cite{Halp62} for 
$m>2$. Strong convergence to $B(r;u_1,\ldots,u_m)$ in the general 
affine case ($u_i\not\equiv 0$) is a routine modification of 
Halperin's proof via Lemma~\ref{l:45} (see \cite{Deut01} for a 
detailed account). Interestingly, if the projectors are not activated 
periodically in \eqref{e:el-nido2009-03-10a} but in a more chaotic
fashion, only weak convergence has been established \cite{Amem65} 
and it is still an open question whether strong convergence holds.

In connection with \eqref{e:el-nido2009-03-10a}, an important question 
is whether the convergence of $(x_n)_{n\in\NN}$ to 
$B(r;u_1,\ldots,u_m)$ occurs at a linear 
rate. The answer is negative and it has actually been shown that 
arbitrarily slow convergence may occur 
\cite{Baus09} in the sense that, for every sequence 
$(\alpha_n)_{n\in\NN}$ in $\zeroun$ such that
$\alpha_n\downarrow 0$, there exits $r\in\HH$ such that
\begin{equation}
\label{e:el-nido2009-03-10c}
(\forall n\in\NN)\quad\|x_n-B(r;u_1,\ldots,u_m)\|\geq\alpha_n.
\end{equation}
On the other hand, several conditions have been found 
\cite{Baus96,Baus09,Baus03,Deut84,Deut08,Kaya88} that guarantee that,
if \eqref{e:puertoprincessa-mai2008-1} admits a solution for some
$(u_i)_{1\leq i\leq m}\in\cart_{i=1}^mU_i$, then, for every $r\in\HH$,
the sequence $(x_n)_{n\in\NN}$ generated by \eqref{e:el-nido2009-03-10a}
converges uniformly linearly to $B(r;u_1,\ldots,u_m)$ 
in the sense that there exists $\alpha\in\Zeroun$ such that 
\cite[Section~4]{Deut08}
\begin{equation}
\label{e:puertoprincessa-mars2009-2}
(\forall n\in\NN)\quad
\|x_n-B(r;u_1,\ldots,u_m)\|\leq\alpha^n\|r-B(r;u_1,\ldots,u_m)\|.
\end{equation}
The next result states that the IBAP implies uniform linear convergence 
of the periodic projection algorithm for solving the underlying affine 
feasibility problem \eqref{e:puertoprincessa-mai2008-1} for every
$(u_i)_{1\leq i\leq m}\in\cart_{i=1}^mU_i$ and every $r\in\HH$. 
In other words, if \eqref{e:puertoprincessa-mai2008-1} admits a 
solution for every $(u_i)_{1\leq i\leq m}\in\cart_{i=1}^mU_i$, then 
uniform linear convergence always occurs in 
\eqref{e:el-nido2009-03-10a}.

\begin{proposition}
\label{p:palawan32}
Suppose that $(U_i)_{1\leq i\leq m}$ satisfies the inverse 
best approximation property and set
\begin{equation}
\label{e:el-nido2009-03-07b}
\alpha=\sqrt{1-\prod_{i=1}^{m-1}\Big(1-{\mathsf c}
\big(U_i^\bot,U_{i+}^\bot\big)^2\Big)}.
\end{equation}
Then $\alpha\in\Zeroun$ and, for every $r\in\HH$ and every
$(u_i)_{1\leq i\leq m}\in\cart_{i=1}^mU_i$, the sequence 
$(x_n)_{n\in\NN}$ generated by \eqref{e:el-nido2009-03-10a} satisfies 
\eqref{e:puertoprincessa-mars2009-2}.
\end{proposition}
\begin{proof}
We first deduce from the implication
\ref{t:palawan-mai2008i}$\Rightarrow$\ref{t:palawan-mai2008vi} 
in Theorem~\ref{t:palawan-mai2008} that
$(\forall i\in\{1,\ldots,m-1\})$ ${\mathsf c}(U_i,U_{i+})<1$.
Hence, it follows from \cite[Theorem~9.35]{Deut01} that
$(\forall i\in\{1,\ldots,m-1\})$ ${\mathsf c}(U^\bot_i,U^\bot_{i+})<1$.
In turn, \eqref{e:el-nido2009-03-07b} and 
\eqref{e:palawan-mai2008} imply that 
\begin{equation}
\label{e:palawan08}
\alpha=\sqrt{1-\prod_{i=1}^{m-1}\bigg(1-{\mathsf c}
\bigg(U_i^\bot,\bigcap_{j=i+1}^mU_j^\bot\bigg)^2\bigg)}\in\Zeroun.
\end{equation}
Now let $(u_i)_{1\leq i\leq m}\in\cart_{i=1}^mU_i$. 
Since the IBAP holds, we have 
\begin{equation}
\label{e:rio-10mai2009}
S(u_1,\ldots,u_m)\neq\emp.
\end{equation}
Altogether, it follows from \eqref{e:palawan08}, 
\eqref{e:rio-10mai2009}, and \cite[Corollary~9.34]{Deut01} applied
to $(U_i^\bot)_{1\leq i\leq m}$ that
\eqref{e:puertoprincessa-mars2009-2} holds.
\end{proof}

In the case when $m=2$, the above result admits a partial converse
based on a result of \cite{Baus09}.

\begin{proposition}
\label{p:diliman2009-03-12}
Suppose that $U_1\cap U_2=\{0\}$, that $(U_1,U_2)$ does not 
satisfy the IBAP, and that $(u_1,u_2)\in U_1\times U_2$ satisfies
$S(u_1,u_2)\neq\emp$. Let
$(\alpha_n)_{n\in\NN}$ be a sequence in $\zeroun$ such that
$\alpha_n\downarrow 0$. Then there exits $r\in\HH$ such that
the sequence $(x_n)_{n\in\NN}$ generated 
by \eqref{e:el-nido2009-03-10a} with $m=2$ satisfies
\begin{equation}
\label{e:el-nido2009-03-10d}
(\forall n\in\NN)\quad\|x_n-B(r;u_1,u_2)\|\geq\alpha_n.
\end{equation}
\end{proposition}
\begin{proof}
It follows from our hypotheses and the equivalence 
\ref{c:palawan-mai2008i}$\Leftrightarrow$\ref{c:palawan-mai2008v} in 
Corollary~\ref{c:palawan-mai2008} that $U_1+U_2$ is not closed.
In turn, we derive from \cite[Theorem~1.4(2)]{Baus09} that there 
exists $y_0\in\HH$ such that the sequence $(y_n)_{n\in\NN}$ 
generated by the alternating projection algorithm
\begin{equation}
\label{e:diliman2009-03-12b}
(\forall n\in\NN)\quad y_{n+1}=P_{U_1^\bot}P_{U_2^\bot}y_n
\end{equation}
satisfies
\begin{equation}
\label{e:diliman2009-03-12c}
(\forall n\in\NN)\quad\|y_n-P_{U_1^\bot\cap U_2^\bot}y_0\|\geq\alpha_n.
\end{equation}
Now let $y\in S(u_1,u_2)$ and set $r=y+y_0$. It follows from 
Proposition~\ref{p:carac}\ref{p:caracii} that 
$S(u_1,u_2)=y+(U_1^\bot\cap U_2^\bot)$. 
Hence, it follows from \eqref{e:diliman2009-03-12a} and
Lemma~\ref{l:45} that
\begin{equation}
\label{e:rio-9mai2009}
B(r;u_1,u_2)=y+P_{U_1^\bot\cap U_2^\bot}(r-y)
=y+P_{U_1^\bot\cap U_2^\bot}y_0. 
\end{equation}
On the other hand, $x_0-y=y_0$ and, 
using Lemma~\ref{l:45}, \eqref{e:el-nido2009-03-10a} with $m=2$ 
and \eqref{e:rio-9mai2009} yield
\begin{equation}
\label{e:diliman2009-03-12d}
(\forall n\in\NN)\quad
x_{n+1}-y=P_{u_1+U_1^\bot}P_{u_2+U_2^\bot}x_n-y
=P_{y+U_1^\bot}P_{y+U_2^\bot}x_n-y
=P_{U_1^\bot}P_{U_2^\bot}(x_n-y).
\end{equation}
This and \eqref{e:diliman2009-03-12b} imply by induction that 
$(\forall n\in\NN)$ $x_n-y=y_n$. In turn, we derive from 
\eqref{e:rio-9mai2009} and \eqref{e:diliman2009-03-12c} that
\begin{equation}
\label{e:diliman2009-03-12e}
(\forall n\in\NN)\quad
\|x_n-B(r;u_1,u_2)\|
=\|(y_n+y)-(y+P_{U_1^\bot\cap U_2^\bot}y_0)\|
=\|y_n-P_{U_1^\bot\cap U_2^\bot}y_0\|
\geq\alpha_n,
\end{equation}
which completes the proof.
\end{proof}

\section{Applications}
\label{sec:4}

In this section, we present several applications of 
Theorem~\ref{t:palawan-mai2008}. As usual, $L^2(\RR^N)$
is the space of real- or complex-valued absolutely
square-integrable functions
on the $N$-dimensional Euclidean space $\RR^N$, 
$\widehat{x}$ denotes the Fourier transform of a function
$x\in L^2(\RR^N)$ and $\supp\widehat{x}$ the support of
$\widehat{x}$. Moreover,
if $A\subset\RR^N$, $1_A$ denotes the characteristic function 
of $A$ and $\complement A$ the complement of $A$. Finally, 
$\mu$ designates the Lebesgue measure on $\RR^N$, $\ran T$ the 
range of an operator $T$ and $\cran T$ is the closure of $\ran T$.

The following lemma and its subsequent refinement, will be used on 
several occasions.

\begin{lemma}{\rm\cite[Proposition~8]{Amre77}, \cite[Corollary~1]{Bene85}}
\label{l:1}
Let $A$ and $B$ be measurable subsets of $\RR^N$ of finite Lebesgue 
measure, and let $x\in L^2(\RR^N)$ be such that $x1_{\complement A}=0$ 
and $\widehat{x}1_{\complement B}=0$. Then $x=0$.
\end{lemma}

\begin{lemma}{\rm\cite[p.~264]{Amre77}, \cite[Theorem~8.4]{Foll97}}
\label{l:2}
Let $A$ and $B$ be measurable subsets of $\RR^N$ of finite Lebesgue
measure. Set $U=\menge{x\in L^2(\RR^N)}{x1_{\complement A}=0}$ and
$V=\menge{x\in L^2(\RR^N)}{\widehat{x}1_{\complement B}=0}$. Then
$\|P_UP_V\|<1$.
\end{lemma}

\subsection{Systems of linear equations}

Going back to Definition~\ref{d:el-nido2009-03-07}, we can say that 
$(U_i)_{1\leq i\leq m}$ satisfies the IBAP if for every 
$(u_i)_{1\leq i\leq m}\in\cart_{i=1}^m\ran P_i$ there exists 
$x\in\HH$ such that $(\forall i\in\{1,\ldots,m\})$ $P_ix=u_i$.
As we have shown, this property holds if \ref{t:palawan-mai2008iv} 
in Theorem~\ref{t:palawan-mai2008} is satisfied, i.e., if 
$(\forall i\in\{1,\ldots,m-1\})$ $\ker P_i+\bigcap_{j=i+1}^m\ker P_j=\HH$.
In the following proposition, we show that such surjectivity results
remain valid if projectors are replaced by more general linear
operators.

\begin{proposition}
\label{p:el-nido2009-03-04}
For every $i\in\{1,\ldots,m\}$, let $\GG_i$ be a normed vector space and 
let $T_i\colon\HH\to\GG_i$ be linear and bounded. Suppose that
\begin{equation}
\label{e:el-nido2009-03-04a}
(\forall i\in\{1,\ldots,m-1\})\quad
\ker T_i+\bigcap_{j=i+1}^m\ker T_j=\HH.
\end{equation}
Then, for every $(y_i)_{1\leq i\leq m}\in\cart_{i=1}^m\ran T_i$,
there exists $x\in\HH$ such that 
\begin{equation}
\label{e:el-nido2009-03-04b}
(\forall i\in\{1,\ldots,m\})\quad T_ix=y_i.
\end{equation}
\end{proposition}
\begin{proof}
For every $i\in\{1,\ldots,m\}$, let
$y_i\in\ran T_i$, set $U_i=(\ker T_i)^\bot$, and let 
$u_i\in U_i$ be such that $T_iu_i=y_i$. Now let
$x\in\HH$. Then $x$ solves \eqref{e:el-nido2009-03-04b} 
$\Leftrightarrow$ 
$(\forall i\in\{1,\ldots,m\})$ $T_ix=T_iu_i$
$\Leftrightarrow$ 
$(\forall i\in\{1,\ldots,m\})$ $T_i(x-u_i)=0$
$\Leftrightarrow$ 
$(\forall i\in\{1,\ldots,m\})$ $x-u_i\in\ker T_i=U_i^\bot$
$\Leftrightarrow$ 
$(\forall i\in\{1,\ldots,m\})$ $P_ix=u_i$.
We thus recover an instance problem \eqref{e:puertoprincessa-mai2008-1} 
and, in view of the equivalence between items \ref{t:palawan-mai2008i} 
and \ref{t:palawan-mai2008iv} in Theorem~\ref{t:palawan-mai2008}, we 
obtain the existence of solutions to \eqref{e:el-nido2009-03-04b} if, 
for every $i\in\{1,\ldots,m-1\}$, $U_i^\bot+U_{i+}^\bot=\HH$, i.e., if
\eqref{e:el-nido2009-03-04a} holds.
\end{proof}

We now give an application of Proposition~\ref{p:el-nido2009-03-04}
to systems of integral equations.

\begin{proposition}
\label{p:el-nido2009-03-05}
For every $i\in\{1,\ldots,m\}$, let $v_i$, $w_i$, and $y_i$ be 
functions in $L^2(\RR^N)$ such that there exists $x_i\in L^2(\RR^N)$
that satisfies $\int_{\RR^N}x_i(s)v_i(s)w_i(t-s)ds=y_i(t)$ 
$\mu$-a.e.\ on \ $\RR^N$. Moreover, suppose that there exist measurable 
sets $(A_i)_{1\leq i\leq m}$ in $\RR^N$ such that  
\begin{equation}
\label{e:supports}
(\forall i\in\{1,\ldots,m\})\quad
\mu\big((A_i+\supp\widehat{v_i})\cap\supp\widehat{w_i}\big)=0 
\end{equation}
and
\begin{equation}
\label{e:overlapping}
(\forall i\in\{1,\ldots,m-1\})\quad
A_i\cup\bigcap_{j=i+1}^m A_j=\RR^N. 
\end{equation}
Then there exists $x\in L^2(\RR)$ such that 
\begin{equation}
\label{e:el-nido2009-03-05}
(\forall i\in\{1,\ldots,m\})\quad
\int_{\RR^N}x(s)v_i(s)w_i(t-s)ds=y_i(t)\;\:\mu\text{-a.e.\ on\;}\RR^N. 
\end{equation}
\end{proposition}
\begin{proof}
The result is an application of 
Proposition~\ref{p:el-nido2009-03-04} in $\HH=L^2(\RR^N)$.
To see this, denote by $\star$ the $N$-dimensional convolution operator 
and, for every $i\in\{1,\ldots,m\}$ and every $x\in\HH$, set 
$T_ix=(xv_i)\star w_i$. Then $(T_i)_{1\leq i\leq m}$ are bounded linear 
operators from $\HH$ to $\HH$ since, by
\cite[Th\'eor\`eme~IV.15]{Brez93},
\begin{equation}
\big(\forall i\in\{1,\ldots,m\}\big)\big(\forall x\in \HH\big)
\quad\|T_ix\|=\|(xv_i)\star w_i\|\leq\|xv_i\|_{L^1}\|w_i\|
\leq\|x\|\,\|v_i\|\,\|w_i\|.
\end{equation}
Now fix $i\in\{1,\ldots,m-1\}$. Since \eqref{e:el-nido2009-03-05} 
can be written as \eqref{e:el-nido2009-03-04b},
Proposition~\ref{p:el-nido2009-03-04} asserts that it suffices to 
show that
\begin{equation}
\label{e:el-nido2009-03-05b}
\ker T_i+\bigcap_{j=i+1}^m\ker T_j=\HH.
\end{equation}
To this end, let $z\in \HH$. It follows from 
\eqref{e:overlapping} that we can write $z=z_1+z_2$, where
$\widehat{z}_1=\widehat{z}\,1_{A_i}$ and 
$\widehat{z}_2=\widehat{z}\,1_{\complement A_i} $. 
We have
\begin{equation}
\label{e:elnido-mars2009-f}
\widehat{T_iz_1}=\big[(z_1v_i)\star w_i\big]^\wedge=
(\widehat{z_1}\star\widehat{v_i})\widehat{w_i}=
\big((\widehat{z}\,1_{A_i})\star\widehat{v_i}\big)\widehat{w_i}
\end{equation}
and
\begin{equation}
\label{e:elnido-mars2009-g}
\supp\big((\widehat{z}\,1_{A_i})\star\widehat{v_i}\big)\subset
\supp(\widehat{z}\,1_{A_i})+\supp\widehat{v_i}\subset
A_i+\supp\widehat{v_i}.
\end{equation}
Therefore, we derive from \eqref{e:elnido-mars2009-f} and 
\eqref{e:supports} that
\begin{equation}
\mu\big(\supp\widehat{T_iz_1}\big)=
\mu\Big(\supp\big((\widehat{z}\,1_{A_i})\star\widehat{v_i}\big)
\cap\supp\widehat{w_i}\Big)=0.
\end{equation}
This shows that $z_1\in\ker T_i$. Now fix $j\in\{i+1,\ldots,m\}$.
Then it remains to show that $z_2\in\ker T_j$. Since 
\eqref{e:overlapping} yields 
$\complement A_i=\bigcap_{k=i+1}^mA_k\subset A_j$,
arguing as above, we get
\begin{equation}
\supp\widehat{T_jz_2}=
\supp\big(((\widehat{z}\,1_{\complement A_i})\star\widehat{v_j})
\widehat{w_j}\big)
\subset\big(\complement
A_i+\supp\widehat{v_j}\big)\cap\supp\widehat{w_j}
\subset\big(A_j+\supp\widehat{v_j}\big)\cap\supp\widehat{w_j}.
\end{equation}
In turn, we deduce from \eqref{e:supports} that
$\mu(\supp\widehat{T_jz_2})=0$ and therefore that
$z_2\in\ker T_j$.
\end{proof} 

We now give an example in which the hypotheses of 
Proposition~\ref{p:el-nido2009-03-05} are satisfied with $m=3$.

\begin{example} 
Let $\{\alpha,\beta,\gamma\}\subset\RR$ and let
$\{v_1,v_2,v_3,w_1,w_2,w_3\}\subset L^{2}(\RR)$. Suppose that
$0<\gamma<2\alpha$ and that
\begin{equation}
\begin{cases}
\supp\widehat{v_1}\subset\left[\beta,\beta+\gamma\right], \;
\supp\widehat{v_2}\subset\left[\alpha,\pinf\right[, \;
\supp\widehat{v_3}\subset\left]\minf,-\alpha\right]\\
\supp\widehat{w_1}\subset\left[-\alpha+\beta+\gamma,\alpha+\beta\right],\;
\supp\widehat{w_2}\subset\left]\minf,0\right],\;
\supp\widehat{w_3}\subset\left[0,\pinf\right[.
\end{cases}
\end{equation}
Now set $A_1=\left]\minf,-\alpha\right]\cup\left[\alpha,\pinf\right[$, 
$A_2=\left[-\alpha,\pinf\right[$, and $A_3=\left]\minf,\alpha\right]$. 
Then \eqref{e:overlapping} is satisfied and, since
$A_1+\supp\widehat{v_1}\subset\left]\minf,-\alpha+\beta+\gamma\right]
\cup\left[\alpha+\beta,\pinf\right[$,
$A_2+\supp\widehat{v_2}\subset\left[0,\infty\right[$, and 
$A_3+\supp\widehat{v_3}\subset\left]\minf,0\right]$, so is
\eqref{e:supports}.
\end{example}

Next, we consider a moment problem with wavelet frames 
\cite{Chri96,Chri01,Daub92}.

\begin{proposition}
Let $\psi$ be a band-limited function in $L^2(\RR)$, say
$\supp\widehat{\psi}\subset[\rho,\rho]$ for some $\rho\in\RPP$.
Suppose that $(\psi_{j,k})_{(j,k)\in\ZZ^2}$, where 
$\psi_{j,k}\colon t\mapsto 2^{j/2}\psi(2^jt-k)$,
is a frame for $L^2(\RR)$, i.e., there exist 
constants $\alpha$ and $\beta$ in $\RPP$ such that
\begin{equation}
\label{e:frame}
\big(\forall x\in L^2(\RR)\big)\quad
\alpha\|x\|^2\leq\sum_{j\in\ZZ}\sum_{k\in\ZZ}
\abscal{x}{\psi_{j,k}}^2\leq\beta\|x\|^2,
\end{equation}
and, moreover, that $(\psi_{j,k})_{(j,k)\in\ZZ^2}$ admits a lower 
Riesz bound $\gamma\in\RPP$, i.e., 
\begin{equation} 
\label{e:lower-riesz}
\big(\forall (c_{j,k})_{(j,k)\in\ZZ^2}\in\ell^2(\ZZ^2)\big)
\quad\sum_{j\in\ZZ}\sum_{k\in\ZZ}|c_{j,k}|^2\leq 
\gamma\Bigg\|\sum_{j\in\ZZ}\sum_{k\in\ZZ}c_{j,k}\psi_{j,k}\Bigg\|^2.
\end{equation}
Let $A$ be a measurable subset of $\RR$ such that $0<\mu(A)<\pinf$,
let $J\in\ZZ$, and set 
\begin{equation}
\label{e:taichung2009-04-19}
\Lambda=\menge{(j,k)\in\ZZ\times\ZZ}{j\leq J}.
\end{equation}
Then, for every function $y\in L^2(A)$ and every sequence 
$(\eta_{j,k})_{(j,k)\in\Lambda}\in\ell^2(\Lambda)$, there exists 
$x\in L^2(\RR)$ such that
\begin{equation} 
\label{e:desired}
x|_A=y\quad\text{and}\quad(\forall (j,k)\in\Lambda)\quad
\scal{x}{\psi_{j,k}}=\eta_{j,k}.
\end{equation}
\end{proposition}
\begin{proof} 
Set $\HH=L^2(\RR)$, $\GG_1=L^2(A)$, and $\GG_2=\ell^2(\Lambda)$, and 
define bounded linear operators 
\begin{equation}
\label{e:taipei2009-04-25}
T_1\colon\HH\to\GG_1\colon x\mapsto x|_A\quad\text{and}\quad
T_2\colon\HH\to\GG_2\colon x\mapsto 
(\scal{x}{\psi_{j,k}})_{(j,k)\in\Lambda}. 
\end{equation}
Then $\ran T_1=\GG_1$ and, on the other hand, it follows from 
\cite[Lemma~2.2(ii)]{Casa02} and \eqref{e:lower-riesz} that 
$\ran T_2=\GG_2$. Hence, in view of \eqref{e:desired}, we 
must show that, for every 
$y_1\in\ran T_1$ and every $y_2\in\ran T_2$, there exists $x\in\HH$ such 
that $T_1x=y_1$ and $T_2x=y_2$. Appealing to 
Proposition~\ref{p:el-nido2009-03-04}, it is enough to show that 
$\ker T_1+\ker T_2=\HH$ or, equivalently, that
\begin{equation}
\label{e:ker-ker}
U_1^\bot+U_2^\bot=\HH,\quad\text{where}\quad
U_1=\cran T_1^*\quad\text{and}\quad U_2=\cran T_2^*.
\end{equation}
Set $U=\menge{x\in L^2(\RR)}{x1_{\complement A}=0}$,
$B=[-2^J\rho,2^J\rho]$, and
$V=\menge{x\in L^2(\RR)}{\widehat{x}1_{\complement B}=0}$.
By Lemma~\ref{l:1}, $U\cap V=\{0\}$ and it therefore follows
from \cite[Lemma~9.5]{Deut01} and Lemma~\ref{l:2} that 
\begin{equation}
\label{e:taipei2009-04-24}
{\mathsf c}(U,V)=\|P_UP_V-P_{U\cap V}\|=\|P_UP_V\|<1. 
\end{equation}
On the other hand, it follows from \eqref{e:taipei2009-04-25} that
$T_1^*\colon\GG_1\to\HH$ satisfies
\begin{equation}
(\forall y\in\GG_1)(\forall t\in\RR)\quad (T_1^*y)(t)=
\begin{cases}
y(t),&\text{if}\;\;t\in A;\\
0,&\text{otherwise}
\end{cases}
\end{equation}
and that
\begin{equation}
T_2^*\colon\GG_2\to\HH\colon(\eta_{j,k})_{(j,k)\in\Lambda}\mapsto
\sum_{(j,k)\in\Lambda}\eta_{j,k}\psi_{j,k}. 
\end{equation}
Since $\bigcup_{(j,k)\in\Lambda}\supp\widehat{\psi_{j,k}}\subset B$, 
we have
\begin{equation}
U_1\subset U\quad\text{and}\quad U_2\subset V.
\end{equation}
Hence, $U_1\cap U_2=\{0\}$ and 
\eqref{e:taipei2009-04-24} yields
\begin{equation}   
\label{e:angles}
{\mathsf c}(U_1,U_2)\leq{\mathsf c}(U,V)<1.  
\end{equation}
In view of the implication
\ref{c:palawan-mai2008vi}$\Rightarrow$\ref{c:palawan-mai2008iv} in 
Corollary~\ref{c:palawan-mai2008}, we conclude that 
\eqref{e:ker-ker} holds.
\end{proof}

\subsection{Subspaces spanned by nearly pairwise bi-orthogonal 
sequences}

The following proposition provides a wide range of applications of 
Theorem~\ref{t:palawan-mai2008} with $m=3$.
\begin{proposition}
\label{p:example}
Let $(u_{1,k})_{k\in\ZZ}$, $(u_{2,k})_{k\in\ZZ}$, and 
$(u_{3,k})_{k\in\ZZ}$ be orthonormal sequences in $\HH$ such that
\begin{equation}
\label{e:bi-ortho2}
(\forall k\in\ZZ)(\forall i\in\{1,2\})(\forall j\in\{i+1,3\})
(\forall l\in\ZZ\smallsetminus\{k\})
\quad u_{i,k}\perp u_{j,l}. 
\end{equation}
Moreover, suppose that
\begin{equation}
\label{e:neat}
\sup_{k\in\ZZ}\sqrt{\abscal{u_{1,k}}{u_{2,k}}}+
\sup_{k\in\ZZ}\sqrt{\abscal{u_{2,k}}{u_{3,k}}}+
\sup_{k\in\ZZ}\sqrt{\abscal{u_{1,k}}{u_{3,k}}}<1.
\end{equation}
Then, for every sequences $(\alpha_{1,k})_{k\in\ZZ}$, 
$(\alpha_{2,k})_{k\in\ZZ}$, and $(\alpha_{3,k})_{k\in\ZZ}$ in 
$\ell^2(\ZZ)$, there exists $x\in\HH$ such that 
\begin{equation}
(\forall k\in\ZZ)\quad
\alpha_{1,k}=\scal{x}{u_{1,k}},\;
\alpha_{2,k}=\scal{x}{u_{2,k}},\;\text{and}\;\;
\alpha_{3,k}=\scal{x}{u_{3,k}}.
\end{equation}
\end{proposition}
\begin{proof} 
For every $i\in\{1,2,3\}$, set $U_i=\spc\{u_{i,k}\}_{k\in\ZZ}$ and
observe that $(\forall x\in\HH)$ 
$P_ix=\sum_{k\in\ZZ}\scal{x}{u_{i,k}}u_{i,k}$. 
Accordingly, we have to show that $(U_1,U_2,U_3)$ satisfies the IBAP.
Using the equivalence 
\ref{t:palawan-mai2008i}$\Leftrightarrow$\ref{t:palawan-mai2008viii} 
in Theorem~\ref{t:palawan-mai2008}, this
amounts to showing that $\|P_1P_{1+}\|<1$ and $\|P_2P_3\|<1$. 

First, let us fix $i\in\{1,2\}$ and $j\in\{i+1,3\}$, and let us show that
\begin{equation} 
\label{e:quick}
\sup_{k\in\ZZ}\abscal{u_{i,k}}{u_{j,k}}^{2}\leq\|P_iP_j\|  
\leq\sup_{k\in\ZZ}\abscal{u_{i,k}}{u_{j,k}}.
\end{equation}
In view of \eqref{e:bi-ortho2}, we have
\begin{equation}
\label{e:hongkong2009-04-16a}
(\forall x\in\HH)\quad
P_iP_jx=\sum_{l\in\ZZ}\scal{x}{u_{j,l}}\scal{u_{j,l}}{u_{i,l}}u_{i,l}.
\end{equation}
Hence, for every $k\in\ZZ$, 
$P_iP_ju_{i,k}=\abscal{u_{i,k}}{u_{j,k}}^{2}u_{i,k}$ and therefore
$\|P_i P_j\|\geq\|P_iP_ju_{i,k}\|=\abscal{u_{i,k}}{u_{j,k}}^{2}$. 
This proves the first inequality in \eqref{e:quick}. On the other hand,
it follows from \eqref{e:hongkong2009-04-16a} that
\begin{equation}
(\forall x\in\HH)\quad
\|P_iP_jx\|^2=\sum_{l\in\ZZ}|\scal{x}{u_{j,l}}\scal{u_{j,l}}{u_{i,l}}|^2
\leq\sup_{l\in\ZZ}\abscal{u_{j,l}}{u_{i,l}}^2\|x\|^2.
\end{equation}
This proves the second inequality in \eqref{e:quick}. 

Since \eqref{e:neat} and \eqref{e:quick} imply that
$\|P_2P_3\|<1$, it remains to show that $\|P_1P_{1+}\|<1$.
We derive from \eqref{e:neat} and \eqref{e:quick} that
\begin{align}
\label{e:assurance}
(\sqrt{\|P_1P_2\|} +\sqrt{\|P_1P_3\|})(1+\sqrt{\|P_2P_3\|})
&=\sqrt{\|P_1 P_2\|}+\sqrt{\|P_1P_3\|} -\|P_2P_3\|   \nonumber\\
&\quad +(\sqrt{\|P_1 P_2\|}+\sqrt{\|P_1P_3\|}
+\sqrt{\|P_2P_3\|})\sqrt{\|P_2P_3\|}\nonumber\\
&<1-\|P_2P_3\|.
\end{align}
For every $k\in\ZZ$, let $P^{\bot}_{3,k}$ denote 
the projector onto $\{u_{3,k}\}^\bot$ and set
\begin{equation}
\label{e:puerto-coffee-mars2009b}
v_{2,k}=\frac{P^\bot_{3,k}u_{2,k}}{\|P^\bot_{3,k}u_{2,k}\|}
=\frac{u_{2,k}-\scal{u_{2,k}}{u_{3,k}}u_{3,k}}
{\sqrt{1-\abscal{u_{2,k}}{u_{3,k}}^2}},
\end{equation}
which is well defined since \eqref{e:neat} guarantees that
$\abscal{u_{2,k}}{u_{3,k}}<1$. Let us note that
\eqref{e:puerto-coffee-mars2009b} yields
\begin{align}
\label{e:2009-04-08}
u_{3,k}-\frac{\scal{u_{3,k}}{u_{2,k}}v_{2,k}}
{\sqrt{1-\abscal{u_{2,k}}{u_{3,k}}^2}}
&=u_{3,k}-\frac{\scal{u_{3,k}}{u_{2,k}}u_{2,k}}
{1-\abscal{u_{2,k}}{u_{3,k}}^2}
+\frac{\abscal{u_{2,k}}{u_{3,k}}^2u_{3,k}}
{1-\abscal{u_{2,k}}{u_{3,k}}^2}\nonumber\\
&=\frac{1}{1-\abscal{u_{2,k}}{u_{3,k}}^2}
\big(u_{3,k}-\scal{u_{3,k}}{u_{2,k}}u_{2,k}\big).
\end{align}
On the other hand, it follows from \eqref{e:bi-ortho2} and 
\eqref{e:puerto-coffee-mars2009b} that 
$\{v_{2,k}\}_{k\in\ZZ}\cup\{u_{3,k}\}_{k\in\ZZ}$ 
is an orthonormal set and that
\begin{equation}
\label{e:puertoprincesa-coffeeshop-mars2009}
\spc\big(\{v_{2,k}\}_{k\in\ZZ}\cup\{u_{3,k}\}_{k\in\ZZ}\big)=
\spc\big(\{P^\bot_{3,k}u_{2,k}\}_{k\in\ZZ}\cup\{u_{3,k}\}_{k\in\ZZ}\big)
=\overline{U_2+U_3}=\overline{U_{1+}}.
\end{equation}
To compute $\|P_1P_{1+}\|$, let $x\in\HH$ and let $k\in\ZZ$. We 
derive from \eqref{e:puertoprincesa-coffeeshop-mars2009} that 
\begin{equation}
P_{1+}x=\sum_{l\in\ZZ}\scal{x}{v_{2,l}}v_{2,l} 
+\sum_{l\in\ZZ}\scal{x}{u_{3,l}}u_{3,l}.
\end{equation}
Hence, using \eqref{e:bi-ortho2}, \eqref{e:puerto-coffee-mars2009b}, 
and \eqref{e:2009-04-08}, we obtain
\begin{align}
\label{e:coeff}
\scal{P_{1+} x}{u_{1,k}}  
&=\scal{x}{v_{2,k}}\scal{v_{2,k}}{u_{1,k}} 
+\scal{x}{u_{3,k}}\scal{u_{3,k}}{u_{1,k}}\nonumber\\
&=\scal{x}{u_{2,k}}\frac{\scal{v_{2,k}}{u_{1,k}}}
{\sqrt{1-\abscal{u_{2,k}}{u_{3,k}}^2}} \nonumber\\
&\quad\;+\scal{x}{u_{3,k}}\Bigg(\scal{u_{3,k}}{u_{1,k}}-
\frac{\scal{u_{3,k}}{u_{2,k}}\scal{v_{2,k}}{u_{1,k}}}
{\sqrt{1-\abscal{u_{2,k}}{u_{3,k}}^2}}\Bigg)\nonumber\\
&=\scal{x}{u_{2,k}}\beta_k+\scal{x}{u_{3,k}}\gamma_k,
\end{align}
where
\begin{equation}  
\beta_k=\frac{\scal{u_{2,k}}{u_{1,k}}-\scal{u_{2,k}}{u_{3,k}}
\scal{u_{3,k}}{u_{1,k}}}{1-\abscal{u_{2,k}}{u_{3,k}}^2}
\end{equation}
and
\begin{equation}  
\gamma_k=\frac{\scal{u_{3,k}}{u_{1,k}}-
\scal{u_{3,k}}{u_{2,k}}\scal{u_{2,k}}{u_{1,k}}}
{1-\abscal{u_{2,k}}{u_{3,k}}^2}.
\end{equation}
We note that \eqref{e:quick} yields
\begin{equation}
\label{e:betak} 
|\beta_k|\leq\frac{\abscal{u_{1,k}}{u_{2,k}}+
\abscal{u_{2,k}}{u_{3,k}}\abscal{u_{1,k}}{u_{3,k}}}
{1-\abscal{u_{2,k}}{u_{3,k}}^2}\leq\frac{\sqrt{\|P_1P_2\|}+
\sqrt{\|P_2P_3\|}\,\sqrt{\|P_1P_3\|}}{1-\|P_2P_3\|}
\end{equation}
and, likewise,
\begin{equation}
\label{e:gammak}
|\gamma_k|\leq\frac{\sqrt{\|P_1P_3\|}+\sqrt{\|P_2P_3\|}\,
\sqrt{\|P_1P_2\|}}{1-\|P_2P_3\|}.
\end{equation}
Thus, we obtain
\begin{align} 
\label{e:parseval}
(\forall x\in\HH)\quad\|P_1P_{1+}x\|
&=\sqrt{\sum_{k\in\ZZ}\abscal{P_{1+}x}{u_{1,k}}^2}\nonumber\\
&\leq\sqrt{\sum_{k\in\ZZ}|\scal{x}{u_{2,k}}\beta_k|^2}
+\sqrt{\sum_{k\in\ZZ}|\scal{x}{u_{3,k}}\gamma_k|^2}\nonumber\\
&\leq\bigg(\sup_{k\in\ZZ}|\beta_k|+\sup_{k\in\ZZ}|\gamma_k|\bigg)
\|x\|\nonumber\\
&\leq\frac{(\sqrt{\|P_1 P_2\|}+\sqrt{\|P_1P_3\|})
(1+\sqrt{\|P_2P_3\|})}{1-\|P_2P_3\|}\|x\|. 
\end{align}
Appealing to \eqref{e:assurance}, we conclude that $\|P_1P_{1+}\|<1$.
\end{proof}

\begin{remark}
A concrete example of subspaces satisfying the hypotheses of
Proposition~\ref{p:example} can be constructed from an orthonormal 
wavelet basis. Take $\psi\in L^2(\RR)$ such that 
the functions $(\psi_{k,l})_{k\in\ZZ^2}$, where
$\psi_{k,l}\colon t\mapsto 2^{k/2}\psi(2^{k}t-l)$, 
form an orthonormal basis of $L^2(\RR)$ \cite{Daub92}. 
For every $i\in\{1,2,3\}$ let, for every $k\in\ZZ$,  
$(\eta_{i,k,l})_{l\in\ZZ}$ be a sequence in $\ell^2(\ZZ)$ such that
$\sum_{l\in\ZZ}|\eta_{i,k,l}|^2=1$ and define
\begin{equation}
U_i=\spc\{u_{i,k}\}_{k\in\ZZ},\quad\text{where}\quad
(\forall k\in\ZZ)\quad u_{i,k}=\sum_{l\in\ZZ}\eta_{i,k,l}\psi_{k,l}.
\end{equation}
Then $(u_{1,k})_{k\in\ZZ}$, $(u_{2,k})_{k\in\ZZ}$, and 
$(u_{3,k})_{k\in\ZZ}$ are orthonormal sequences in $L^2(\RR)$ that
satisfy \eqref{e:bi-ortho2}. Moreover since, for every $i$ and $j$ 
in $\{1,2,3\}$ and every $k\in\ZZ$, $\scal{u_{i,k}}{u_{j,k}}=
\sum_{l\in\ZZ}\eta_{i,k,l}\overline{\eta_{j,k,l}}$,
the main hypothesis \eqref{e:neat} is equivalent to
\begin{equation}
\sup_{k\in\ZZ}\sqrt{\bigg|\sum_{l\in\ZZ}\eta_{1,k,l}
\overline{\eta_{2,k,l}}\bigg|}+
\sup_{k\in\ZZ}\sqrt{\bigg|\sum_{l\in\ZZ}
\eta_{2,k,l}\overline{\eta_{3,k,l}}\bigg|}+
\sup_{k\in\ZZ}\sqrt{\bigg|\sum_{l\in\ZZ}\eta_{1,k,l}
\overline{\eta_{3,k,l}}\bigg|}<1.
\end{equation}
\end{remark}

\subsection{Harmonic analysis and signal recovery}

Many problems arising in areas such as harmonic analysis
\cite{Amre77,Bene85,Foll97,Havi94,Jami07,Mele96}, signal theory
\cite{Byrn05,Papo75,Youl78}, image processing \cite{Proc93,Star87},
and optics \cite{Mont82,Star81} involve imposing known values of an
ideal function in the time (or spatial) and Fourier domains. 
In this section, we describe applications of 
Theorem~\ref{t:palawan-mai2008} to such problems.

The following lemma will be required.

\begin{lemma}
\label{l:3}
Let $U$, $V$, and $W$ be closed vector subspaces of $\HH$ such that
$W\subset V$. Then $\|P_UP_W\|\leq\|P_UP_V\|$.
\end{lemma}
\begin{proof}
Set $B=\menge{x\in\HH}{\|x\|\leq 1}$. Then $P_W(B)\subset B$. In turn, 
since $W\subset V$, $P_W(B)=P_V(P_W(B))\subset P_V(B)$ and hence
$P_U(P_W(B))\subset P_U(P_V(B))$. Consequently,
$\|P_UP_W\|=\sup \menge{\|P_U P_W x \|}{x\in B}        
\leq\sup\menge{\|P_U P_V x\|}{x\in B}
=\|P_UP_V\|$.
\end{proof}

The scenario of the next proposition has a simple interpretation in
signal recovery \cite{Proc93,Star87}: an $N$-dimensional square-summable 
signal has known values over certain domains of the spatial and 
frequency domains and, in addition, $m-2$ scalar linear measurements of 
it are available.

\begin{proposition}
\label{p:77}
Let $A$ and $B$ be measurable subsets of $\RR^N$ of finite Lebesgue
measure, and suppose that $m\geq 3$. Moreover, let 
$(v_i)_{1\leq i\leq m-2}$ be functions in $L^2(\RR^N)$ with 
disjoint supports $(C_i)_{1\leq i\leq m-2}$ such that    
\begin{equation}
\label{e:77}
(\forall i\in\{1,\ldots,m-2 \})\quad\mu(C_i)<\pinf\quad\text{and}\quad
\mu(C_i\cap\complement A)>0.
\end{equation}
Then, for every functions $v_m$ and $v_{m-1}$ in $L^2(\RR^N)$ and every 
$(\eta_i)_{1\leq i\leq m-2}\in\RR^{m-2}$, there exists a function 
$x\in L^2(\RR^N)$ such that 
\begin{equation}
\label{e:rainy77}
(\forall i\in\{1,\ldots,m-2 \})\quad
\int_{C_i}x(t)\overline{v_i(t)}dt=\eta_i,\;\;
x|_A=v_{m-1}|_A,\;\;\text{and}\;\;
\widehat{x}|_{B}=\widehat{v}_m|_{B}.
\end{equation}
\end{proposition}
\begin{proof}
We first observe that the problem under consideration is a special case of 
\eqref{e:puertoprincessa-mai2008-1} with $\HH=L^2(\RR^N)$,
\begin{equation}
\label{e:kokuyo1}
\begin{cases}
U_i=\spa\{v_i\}&\text{and}\quad
u_i=\eta_iv_i/\|v_i\|^2,\quad 1\leq i\leq m-2; \\
U_{m-1}=\menge{x\in\HH}{x1_{\complement A}=0}
&\text{and}\quad u_{m-1}=v_{m-1}1_A;  \\
U_m=\menge{x\in\HH}{\widehat{x}1_{\complement B}=0}
&\text{and}\quad\widehat{u_m}=\widehat{v_m}1_B.
\end{cases}
\end{equation}
It follows from Lemma~\ref{l:2} that $\|P_{m-1}P_m\|<1$. Hence, in view 
of Corollary~\ref{c:paris-octobre2008}, it suffices to show that the 
closed vector subspaces $(U_i)_{1\leq i\leq m}$ are linearly 
independent. Since the supports 
$(C_i)_{1\leq i\leq m-2}$ are disjoint, the subspaces
$(U_i)_{1\leq i\leq m-2}$ are independent. Therefore, if we set
$U=\sum_{i=1}^{m-2}U_i$, it is enough to show that $U$, $U_{m-1}$, 
and $U_m$ are independent. To this end, take 
$(y,y_{m-1},y_m)\in U\times U_{m-1}\times U_m$ such that 
\begin{equation}
\label{e:anr05}
y+y_{m-1}+y_m=0,
\end{equation}
and set $C=\bigcup_{i=1}^{m-2}C_i$. We have
$(y+y_{m-1})1_{\complement(A\cup C)}=0$, $\mu(A\cup C)<\pinf$,
$\widehat{y}_m1_{\complement B}=0$, and $\mu(B)<\pinf$. 
Hence, it follows from \eqref{e:anr05} and Lemma~\ref{l:1} that 
\begin{equation}
\label{e:kokuyo2}
y+y_{m-1}=0\;\;\text{and}\;\;y_m=0.
\end{equation}
It remains to show that $y=0$. Since $y\in U$, there exist
$(\alpha_i)_{1\leq i\leq m-2}\in\CC^{m-2}$ such that 
$y=\sum_{i=1}^{m-2}\alpha_iv_i$. However, since the supports
$(C_i)_{1\leq i\leq m-2}$ are disjoint,
\begin{equation}
\label{e:kokuyo5}
\|y\|^2=\bigg\|\sum_{i=1}^{m-2}\alpha_iv_i\bigg\|^2
=\sum_{i=1}^{m-2}|\alpha_i|^2\|v_i\|^2.
\end{equation}
On the other hand, \eqref{e:77} implies that, for every 
$i\in\{1,\ldots,m-2\}$,
\begin{equation}
\label{e:kokuyo4}
\|v_i\|^2=\int_{C_i\cap A}|v_i(t)|^2dt
+\int_{C_i\cap\complement A}|v_i(t)|^2dt
>\int_{C_i\cap A}|v_i(t)|^2dt=\|v_i1_A\|^2.
\end{equation}
At the same time, we derive from 
\eqref{e:kokuyo2} that $y=-y_{m-1}\in U_{m-1}$ and therefore from 
\eqref{e:kokuyo1} that $y1_{\complement A}=0$. Consequently, 
\eqref{e:kokuyo5} yields
\begin{equation}
\label{e:kokuyo3}
\sum_{i=1}^{m-2}|\alpha_i|^2\|v_i\|^2=
\|y\|^2=\|y1_A\|^2=\bigg\|\sum_{i=1}^{m-2}\alpha_iv_i1_A\bigg\|^2
=\sum_{i=1}^{m-2}|\alpha_i|^2\|v_i1_A\|^2.
\end{equation}
In view of \eqref{e:kokuyo4}, we conclude that
$(\forall i\in\{1,\ldots,m-2\})$ $\alpha_i=0$.
\end{proof}

\begin{remark}
\label{r:3}
In connection with Proposition~\ref{p:77}, let us make 
a few comments on the following classical problem: given 
measurable subsets $A$ and $B$ of $\RR^N$ such that 
$\mu(A)>0$ and $\mu(B)>0$, and functions $a$ and $b$
in $L^2(\RR^N)$, is there a function $x\in L^2(\RR^N)$ such that
\begin{equation}
\label{e:74-7}
x|_A=a|_A\quad\text{and}\quad\widehat{x}|_{B}=b|_{B}\;?
\end{equation}
To answer this question, let us set
\begin{equation}
\label{e:kokuyo7}
\begin{cases}
U_1=\menge{x\in L^2(\RR^N)}{x1_{\complement A}=0}
&\text{and}\quad u_1=a1_A,\\
U_2=\menge{x\in L^2(\RR^N)}{\widehat{x}1_{\complement B}=0}
&\text{and}\quad\widehat{u_2}=b1_B.
\end{cases}
\end{equation}
Thus, the problem reduces to an instance of
\eqref{e:puertoprincessa-mai2008-1} in which $m=2$.
\begin{itemize}
\item
If $\mu(A)<\pinf$ and $\mu(B)<\pinf$, it follows from Lemma~\ref{l:2}, 
\eqref{e:kokuyo7}, and the implication 
\ref{c:palawan-mai2008viii}$\Rightarrow$\ref{c:palawan-mai2008i} in 
Corollary~\ref{c:palawan-mai2008} that the answer is affirmative
(see also \cite[Corollary~5.B~p.~100]{Havi94}).
\item
If $\mu(\complement A)<\pinf$ and $\mu(\complement B)<\pinf$, it follows
from \eqref{e:kokuyo7}, Proposition~\ref{p:23}, and Lemma~\ref{l:1}
(applied to $U_1^\bot$ and $U_2^\bot$) that \eqref{e:74-7} has at 
most one solution.
\item
Suppose that $A$ is bounded and that $\mu(\complement B)>0$, and 
let $\varepsilon\in\RPP$. Then there exists $x\in L^2(\RR^N)$ such that
\begin{equation}
\label{e:desired-ineq}
\int_A|x(t)-a(t)|^2dt+
\int_B|\widehat{x}(\xi)-b(\xi)|^2d\xi
<\varepsilon.
\end{equation} 
To show this, we first observe that $U_1\cap U_2=\{0\}$. Indeed, let 
$y\in U_1\cap U_2$. Then $\widehat{y}$ can be extended to an entire 
function on $\CC^N$ (see \cite[Theorem~7.23]{Rudi91} or 
\cite[Theorem~III.4.9]{Stei71}) and, at the same time, 
$\widehat{y}1_{\complement B}=0$, which implies that 
$\widehat{y}=0$ \cite[Theorem~I.3.7]{Rang86}. Hence, 
applying Proposition~\ref{p:denseness} with $m=2$, we obtain the 
existence of $x\in L^2(\RR^N)$ such that
\begin{equation} 
\|P_1x-u_1\|^2+\|P_2x-u_2\|^2<\varepsilon, 
\end{equation}
which yields \eqref{e:desired-ineq}. In the case when $\complement B$ 
is a ball centered at the origin and $b=0$, 
\eqref{e:desired-ineq} provides the
following approximate band-limited extrapolation result: there exists
$x\in L^2(\RR^N)$ which approximates $a$ on $A$ and such that
$\widehat{x}$ nearly vanishes for high frequencies.
\end{itemize}
\end{remark}

The following example describes a situation in which the IBAP fails.

\begin{example}
The following example is from \cite{Mont82}. Let 
$C=[-1/2,1/2] \times [-1/2,1/2]$ and set $\HH=L^2(C)$. 
Moreover, define 
\begin{equation}
(\forall (m,n)\in\ZZ^2)\quad\widehat{x}(m,n)=
\int_Cx(s,t)\exp(-i2\pi(ms+nt))dsdt,
\end{equation}
set $A=[0,1/2]\times [0,1/2]$, and set $B=F\cup\menge{(m,0)}{m\in\ZZ}$,
where $F$ is a nonempty finite subset of $\ZZ\times\ZZ$. The problem
amounts to finding functions with prescribed best approximations from 
the closed vector subspaces
\begin{equation}
\begin{cases}
U_1=\menge{x\in\HH}{x1_{\complement A}=0}\\
U_2=\menge{x\in\HH}{x(s,t)=x(-s,t)\;\text{a.e.\ on}\;C}\\
U_3=\menge{x\in\HH}{\widehat{x}1_{B}=0}. 
\end{cases}
\end{equation}
Since $U^{\bot}_1+U^{\bot}_2+U^{\bot}_3=\{0\}$ \cite{Mont82}, it
follows from Proposition~\ref{p:23} that the problem has at most one
solution. However, the subspaces are not independent. Indeed, given 
a finite subset $I$ of $\ZZ$ such that $(0,n)\notin F$ whenever $n\in I$
and complex numbers $(c_n)_{n\in I}$, the trigonometric polynomial
\begin{equation}
(s,t)\mapsto\sum_{n\in I}c_n e^{i2\pi nt} 
\end{equation}
is in $U_2\cap U_3$. Therefore, in the light of 
Corollary~\ref{c:kimono-ken}, the IBAP does not hold.
\end{example}

{\bfseries Acknowledgement:} 
The work of the first author was supported by the Agence Nationale 
de la Recherche under grant ANR-08-BLAN-0294-02.
The work of the second author was supported by the Creative and Research
Scholarship Program of the University of the Philippines System.

\end{document}